\documentclass[10pt,twoside]{article}

\usepackage[centertags]{amsmath}
\usepackage{amsfonts}
\usepackage{amssymb}
\usepackage{amsthm}
\usepackage{epsfig}
\usepackage{color}
\allowdisplaybreaks[4]

\newcommand{\COM}[1]{}

\begin{document}
\baselineskip 15pt \setcounter{page}{1}
\title{\bf \Large  On Piterbarg's max-discretisation theorem for homogeneous Gaussian random
fields\thanks{Research supported by National Science Foundation of
China (No. 11326175), Natural Science Foundation of Zhejiang Province of China (No. LQ14A010012) and Research Start-up Foundation of Jiaxing
University (No. 70512021).} }
\author{{\small Zhongquan Tan}\footnote{ E-mail address:  tzq728@163.com }\\
{\small\it  College of Mathematics, Physics and Information Engineering, Jiaxing University, Jiaxing 314001, PR China}\\
{\small Kaiyong Wang}\\
{\small\it  School of Mathematics and Physics, Suzhou University of Science and Technology, Suzhou 215009, PR China}\\
}
 \maketitle
 \baselineskip 15pt

\begin{quote}
{\bf Abstract:}\ \ Motivated by the papers of Piterbarg (2004) and H\"{u}sler (2004),
in this paper the asymptotic relation between the maximum of a
continuous dependent homogeneous Gaussian random field
and the maximum of this field sampled at discrete time points is
studied. It is shown that, for the weakly dependent case, these two maxima are
asymptotically independent, dependent and coincide when the grid of the discrete time
points is a sparse grid, Pickands grid and dense grid, respectively, while for the strongly dependent case,
 these two maxima are
asymptotically totally dependent if the grid of the discrete time
points is sufficiently dense, and asymptotically dependent if the
the grid points are sparse or Pickands grids.



\end{quote}

\section{Introduction}

Let $\{X(t), t\geq0\}$ be a stationary Gaussian process with mean 0,
variance 1, correlation function $r(t)$ and continuous sample
functions.  The study on the limit distribution theory on the maximum of $\{X(t),
t\geq0\}$ up to time $T$: $M_{T}= \max\{X(t), 0 \leq t \leq T \}$ can be dated back to Pickands (1969).
 Assume that the correlation function $r(t)$
satisfies for some $\alpha\in(0, 2]$,
\begin{eqnarray}
\label{eq1.1}
r(t)=1-|t|^{\alpha}+o(|t|^{\alpha})\ \ \mbox{as}\ \ t\rightarrow 0\ \ \mbox{and}\ \ r(t)<1\ \ \mbox{for}\ \ t>0
\end{eqnarray}
and
\begin{eqnarray}
\label{eq1.2}
r(t)\log t\rightarrow 0,\ \ \mbox{as}\ \  t\rightarrow \infty.
\end{eqnarray}
It is well known (see e.g. Pickands (1969), Leadbetter et al. 1983) that
(\ref{eq1.1}) and (\ref{eq1.2}) imply the following classic limit relation
\begin{eqnarray}
\label{eq1.A}
P\{a_{T}(M_{T}-b_{T})\leq x\}\rightarrow \exp(-e^{-x})
\end{eqnarray}
as $T\rightarrow\infty$, where
$$a_{T}=\sqrt{2\log T}, \ \ b_{T}=\sqrt{2\log T}+\frac{\log[(2\pi)^{-1/2}\mathcal{H}_{\alpha}(2\log T)^{-1/2+1/\alpha}]}{\sqrt{2\log T}}.
$$
Here $\mathcal{H}_{\alpha}$ denotes Pickands constant, which is
defined by
$\mathcal{H}_{\alpha}=\lim_{\lambda\rightarrow\infty}\mathcal{H}_{\alpha}(\lambda)/\lambda$,
with
$$\mathcal{H}_{\alpha}(\lambda)=\mathbb{E}\exp\left(\max_{t\in[0,\lambda]}\sqrt{2}B_{\alpha/2}(t)-t^{\alpha}\right)$$
and $B_{H}$ is a fractional Brownian motion, that is a Gaussian zero
mean process with stationary increments such that
$\mathbb{E}B_{H}^{2}(t)=|t|^{2H}$. It is also well known that
$0<\mathcal{H}_{\alpha}<\infty$, see e.g. Pickands (1969),
Leadbetter et al. (1983), Piterbarg (1996).

The extensions of the classic result (\ref{eq1.A}) to more general cases, such as for
 non-stationary case, strongly dependent case, can be found in
 Mittal and Ylvisaker (1975),  McCormick (1980), McCormick and Qi (2000),
 H\"{u}lser (1990), Konstant and  Piterbarg (1993),
Seleznjev (1991, 1996), H\"{u}lser (1999), H\"{u}lser et al. (2003), Tan et al. (2012) and among others.

In applied fields, however, the classic result (\ref{eq1.A})
can not be used directly, since the available samples are discrete.
Usually, simulation techniques are applied to derive results for
continuous process when they can not be derived with mathematical
analytic tools. Simulations of such processes are performed for
discrete time-grids, while the results should be interpreted in the
context of continuous time. Therefore, it is crucial
to investigate the relation between  the
extremes of the continuous process and the extremes of the discrete
process.

Piterbarg (2004) first studied the asymptotic relation between
$M_{T}$ and the maximum of the discrete version $M_{T}^{\delta}=
\max\{X(k\delta), 0 \leq k\delta \leq T \}$ for some
$\delta=\delta(T)>0, k\in \mathbb{N}$, where $\mathbb{N}$ denotes the set of all natural numbers. Following Piterbarg (2004),
we consider uniform grids
$\mathfrak{R}=\mathfrak{R}(\delta)=\{k\delta: k\in \mathbb{N}\}$,
$\delta>0$. A grid is called sparse if $\delta$ is such that
$$\delta(2\log T)^{1/\alpha}\rightarrow D$$
with $D=\infty$. If $D\in(0, \infty)$, the grid is a Pickands grid, and if $D = 0$, the grid is dense.

For the stationary Gaussian processes, Piterbarg (2004)  first showed that the maximum
$M_{T}^{\delta}$ of discrete time points and the maximum $M_{T}$ of
the continuous time points can be asymptotically independent,
dependent or totally dependent if the grid is a sparse, a Pickands
or a dense grid, respectively. This type of results are called
Piterbarg's max-discretisation theorems in the literature, see eg.  Tan and Hashorva (2014a).

Based on the results of H\"{u}sler (1990), Piterbarg's max-discretisation theorems were
extended by H\"{u}sler (2004) to a class of locally stationary
Gaussian processes which was introduced by Berman (1974). Other
related results such as for the storage process with fractional
Brownian motion as input and stationary non-Gaussian case can be
found in H\"{u}sler and Piterbarg (2004) and Turkman (2012),
respectively. The recent contributions  Tan
and Wang (2013) and Tan and Tang (2014) present Piterbarg's max-discretisation theorem for
strongly dependent stationary  Gaussian processes.
The Piterbarg's max-discretisation theorems for multivariate Gaussian processes
 can be found in Tan and Hashorva (2014b) and their improvement to different grids can be found in
Tan and Hashorva (2015).

The Piterbarg's max-discretisation theorems for Gaussian processes
have been  studied extensively under different conditions in the past, but it is far from complete.
In this paper, we are interested in the similar
problems for the Gaussian random fields. It
is well known that Gaussian random fields play a very important role in many applied sciences, such as in image analysis,
atmospheric sciences, geostatistics, neuroimaging,
astrophysics, oceanography, hydrology and agriculture, among others, see eg. Adler and Taylor (2007)
for details. Extremes and their limit properties are particularly
important in these applications, see, for instance, Aza\"{\i}s and Wschebor (2009).

The paper is organized as follows: In
Section 2, we present the main results for weakly and strongly dependent
Gaussian fields. Section 3 gives the proofs. Some technical auxiliary results are presented in Sections
4 and 5. Let $\phi$ and $\Psi$ denote the density
function and tail distribution function of a standard normal
variable.

\section{Main results}
Denote the set of all real numbers by $\mathbb{R}$ and
let $\mathbb{R}^{d}$ be d-dimensions product space of $\mathbb{R}$,
where $d\geq 2$. In this paper, we only consider the case of $d=2$
since it is notationally simplest and the results for higher
dimensions follow analogouss arguments.
Here the operations with the vectors are meant component-wise. For instance for two vectors
$\mathbf{t}=(t_{1},t_{2})$ and $\mathbf{s}=(s_{1},s_{2})$,
$\mathbf{s}\leq \mathbf{t}$, $\mathbf{t-s}$ and $\mathbf{st}$ mean $s_{i}\leq
t_{i}$, $i=1,2$, $(t_{1}-s_{1},t_{2}-s_{2} )$ and $(s_{1}t_{1},s_{2}t_{2})$, respectively.
 $\mathbf{T}\rightarrow\infty$ means $T_{i}\rightarrow\infty$, $i=1,2$.
Let $\mathbf{I_{T}}=\{\mathbf{t}\in \mathbb{R}^{2}:
0\leq t_{i}\leq T_{i},i=1,2\}$. Let $\{X(\mathbf{t}):\mathbf{t}\geq
\mathbf{0}\}$ denote a homogeneous Gaussian field with covariance
function
$$r(\mathbf{t})=\mathbb{C}ov(X(\mathbf{t}),X(\mathbf{0})).$$ In
this paper we assume that the covariance function satisfies the
following conditions:
\begin{itemize}
\item[\textbf{A1:}]
$r(\mathbf{t})=1-|t_{1}|^{\alpha_{1}}-|t_{2}|^{\alpha_{2}}+o(|t_{1}|^{\alpha_{1}}+|t_{2}|^{\alpha_{2}})$
as $\mathbf{t}\to 0$ with $\alpha_{i}\in(0,2]$;
\item[\textbf{A2:}]
$r(\mathbf{t})<1$ for $\mathbf{t}\neq \mathbf{0}$;
\item[\textbf{A3:}]
$\lim_{\mathbf{t}\rightarrow\infty}r(\mathbf{t})\log(t_{1}t_{2})=r\in[0,\infty)$ and both $r(0,t_{2})\log t_{2}$ and  $r(t_{1},0)\log t_{1}$
are bounded.
%\item[\textbf{A3:}]
%$\sup_{\mathbf{t}\in S(d)}\left|r(\mathbf{t})\log d-r\right|\to
%0$ as $d\to\infty$, where $r\in[0,\infty)$, $S(d)$ denotes the
%sphere of center $(0,0)$ and radius $d>0$ in $\mathbb{R}^2$ with
%Euclidean metric.
\end{itemize}

Throughout the paper, for any set $\mathbf{E}\subset \mathbb{R}^{2}$ and $\mathbf{k}\in \mathbb{N}^{2}$, define
$$M_{\mathbf{E}}=\max_{\mathbf{t}\in \mathbf{E}}X(\mathbf{t})=\max\{X(\mathbf{t}), \mathbf{t}\in \mathbf{E}\},\ \
\ M_{\mathbf{E}}^{\mathbf{p}}=\max_{\mathbf{t}\in
\mathbf{E} \cap \mathfrak{R}(p_{1})\times\mathfrak{R}(p_{2})}X(\mathbf{t})=\max\{X(\mathbf{kp}), \mathbf{kp}\in \mathbf{E}\},$$ where
$\mathfrak{R}(p_{i})=\{kp_{i}, k\in \mathbb{N}\}$, $i=1,2$, are uniform grids.
If $\mathbf{E}=\mathbf{I_{T}}$, we write the above two maxima for simplicity as $M_{\mathbf{T}}$
and $M_{\mathbf{T}}^{\mathbf{p}}$, respectively.
For dealing with the multivariate case, we redefine the uniform grids $\mathfrak{R}(p_{i})=\{kp_{i}, k\in \mathbb{N}\}$, $i=1,2$ as following.
The grid $\mathfrak{R}(p_{i})$ is called sparse if $p_{i}=p_{i}(\mathbf{T})$
is such that
$$p_{i}(2\log T_{1}T_{2})^{1/\alpha_{i}}\rightarrow D_{i},\ \ i=1,2$$
with $D_{i}=\infty$. If $D_{i}\in(0, \infty)$, the grid is a
Pickands grid, and if $D_{i} = 0$, the grid is dense.

Under conditions $\textbf{A1}$ and $\textbf{A2}$, Theorem 7.1 of Piterbarg (1996) showed that for any
fixed $\mathbf{h}>\mathbf{0}$
\begin{eqnarray}
\label{eq2.0}
P\left(\max_{\mathbf{t}\in\mathbf{I_{h}}}X(\mathbf{t})>u\right)
=\mathcal{H}_{\alpha_{1}}\mathcal{H}_{\alpha_{2}}h_{1}h_{2}
u^{2/\alpha_{1}+2/\alpha_{2}}\Psi(u)(1+o(1)),
\end{eqnarray} as
$u\rightarrow\infty$, where $\mathcal{H}_{\alpha_{i}}$, $i=1,2$ are
the Pickands constants. This exact asymptotic plays crucial role in
deriving the Gumbel law and also will be used in the proofs of our
main results.

Now, we state our main results which extend the existing results (including Piterbarg (2004) and Tan and Wang (2013)) to Gaussian random fields.
The extensions are, however, nontrivial in that asymptotic relation between two Gaussian fields
is more complicated than that of Gaussian processes. This can be seen from the proof
that follows.

\textbf{Theorem 2.1}. {\sl Let $\{X(\mathbf{t}):\mathbf{t}\geq
\mathbf{0}\}$ be a homogeneous  Gaussian field with covariance
function $r(\mathbf{t})$ satisfying \textbf{A1}, \textbf{A2} and
\textbf{A3}. Then for any sparse grids $\mathfrak{R}(p_{i})$,
$i=1,2$,
\begin{eqnarray}
\label{eq2.1}
&&P\left\{a_{\mathbf{T}}\big(M_{\mathbf{T}}-b_{\mathbf{T}}\big)\leq x,
           a_{\mathbf{T}}\big(M_{\mathbf{T}}^{\mathbf{p}}- b_{\mathbf{T}}^{\mathbf{p}}\big)\leq y\right\}\nonumber \\
&&\ \ \ \ \ \ \longrightarrow\int_{-\infty}^{+\infty}\exp\left(-\big(e^{-x-r+\sqrt{2r}z}+e^{-y-r+\sqrt{2r}z}\big)\right)\phi(z)dz
\end{eqnarray}
as $\mathbf{T}\rightarrow\infty,$ where
$$a_{\mathbf{T}}=\sqrt{2\log T_{1} T_{2}},\ \
b_{\mathbf{T}}=a_{\mathbf{T}}+a_{\mathbf{T}}^{-1}\log\left((2\pi)^{-1/2}\mathcal{H}_{\alpha_{1}}\mathcal{H}_{\alpha_{2}}
(a_{\mathbf{T}})^{1/\alpha_{1}+1/\alpha_{2}-1/2}\right)$$ and
$$b_{\mathbf{T}}^{\mathbf{p}}=a_{\mathbf{T}}+a_{\mathbf{T}}^{-1}\log\left((2\pi)^{-1/2}p_{1}^{-1}p_{2}^{-1}
(a_{\mathbf{T}})^{-1/2}\right).$$ }

As the special case, we can obtain the limit distribution of the
maximum for  a homogeneous  Gaussian random field, which will be used in the proof of Theorem 2.3.

\textbf{Corollary 2.1}. {\sl Let $\{X(\mathbf{t}):\mathbf{t}\geq
\mathbf{0}\}$ be a homogeneous  Gaussian field with covariance
function $r(\mathbf{t})$ satisfying \textbf{A1}, \textbf{A2} and
\textbf{A3}. Then for any $x\in \mathbb{R}$,
\begin{eqnarray}
\label{eq2.4}
&&P\bigg\{a_{\mathbf{T}}\big(M_{\mathbf{T}}-b_{\mathbf{T}}\big)\leq x\bigg\}
\longrightarrow \int_{-\infty}^{+\infty}\exp\left(-e^{-x-r+\sqrt{2r}z}\right)\phi(z)dz
\end{eqnarray}
as $\mathbf{T}\rightarrow\infty.$ }

\textbf{Remark 2.2}. A similar result as Corollary 2.1 can also be derived from Theorem 15.2 in Chapter 4 of Piterbarg (1996), where the author
dealt with the limit properties of uncrossing point processes under some slight different conditions.

Before presenting the result for Pickands grids, we introduce the
following Pickands type constants.
For $a>0$, define,
$$\mathcal{H}_{a,\alpha}(\lambda)=\mathbb{E}\exp\left(\max_{ka\in[0,\lambda]}\sqrt{2}B_{\alpha/2}(ka)-(ka)^{\alpha}\right),$$
we have (see Leadbetter et al. 1983),
$$\mathcal{H}_{a,\alpha}=\lim_{\lambda\rightarrow\infty}\frac{\mathcal{H}_{a,\alpha}(\lambda)}{\lambda}\in(0,\infty).$$
%Define,
%\begin{eqnarray*}
%&&\mathcal{H}_{a,\alpha}^{x,y}(\lambda)=\int_{-\infty}^{+\infty}e^{s}P\left(\max_{t\in[0,\lambda]}\sqrt{2}B_{\alpha/2}(t)-t^{\alpha}>s+x,
%\max_{k:ka\in [0,\lambda]}\sqrt{2}B_{\alpha/2}(ka)-(ka)^{\alpha}>s+y\right)ds.
%\end{eqnarray*}
For any $\mathbf{d}>\mathbf{0}$, define
\begin{eqnarray*}
&&\mathcal{H}_{\mathbf{d},\alpha_{1},\alpha_{2}}^{x,y}(\lambda_{1},\lambda_{2})=\int_{-\infty}^{+\infty}e^{s}
P\left(\max_{(t_{1},t_{2})\in[0,\lambda_{1}]\times[0,\lambda_{2}]}\sqrt{2}\chi(t_{1},t_{2})>s+x,\right.\\
&&\ \ \ \ \ \ \ \ \ \ \ \ \ \ \ \ \ \ \ \ \ \ \ \ \ \ \ \ \ \ \ \ \ \ \ \ \ \ \ \ \ \
\left.\max_{(k_{1}d_{1},k_{2}d_{2})\in[0,\lambda_{1}]\times[0,\lambda_{2}]}\sqrt{2}\chi(k_{1}d_{1},k_{2}d_{2})>s+y\right)ds,
\end{eqnarray*}
where
$$\chi(t_{1},t_{2})=B^{(1)}_{\alpha_{1}/2}(t_{1}) +
B^{(2)}_{\alpha_{2}/2}(t_{2})-|t_{1}|^{\alpha_{1}}-|t_{2}|^{\alpha_{2}}$$ and
$B^{(1)}_{\alpha_{1}/2}(\cdot)$, $B^{(2)}_{\alpha_{2}/2}(\cdot)$ are two
independent fractional Brownian motions.

\textbf{Theorem 2.2}. {\sl  Let $\{X(\mathbf{t}):\mathbf{t}\geq
\mathbf{0}\}$ be a homogeneous  Gaussian field with covariance
function $r(\mathbf{t})$ satisfying \textbf{A1}, \textbf{A2} and
\textbf{A3}. Then for any Pickands grids $\mathfrak{R}(p_{i})=\mathfrak{R}(a_{i}(2\log
T_{1}T_{2})^{-1/\alpha_{i}})$ with $a_{i}>0$, $i=1,2$, the following limit
exists,
$$\mathcal{H}_{\mathbf{a},\alpha_{1},\alpha_{2}}^{x,y}:=\lim_{\lambda_{1}\rightarrow\infty\atop\lambda_{2}\rightarrow\infty}
\mathcal{H}_{\mathbf{a},\alpha_{1},\alpha_{2}}^{x,y}(\lambda_{1},\lambda_{2})/\lambda_{1}\lambda_{2}\in (0,\infty)$$
and
\begin{eqnarray}
\label{eq2.2}
&&P\left\{a_{\mathbf{T}}\big(M_{\mathbf{T}}-b_{\mathbf{T}}\big)\leq x,
      a_{\mathbf{T}}\big(M_{\mathbf{T}}^{\mathbf{p}}- b_{\mathbf{a},\mathbf{T}}\big)\leq y\right\}\nonumber\\
&&\ \ \ \ \ \longrightarrow \int_{-\infty}^{+\infty}\exp\left(-\big(e^{-x-r+\sqrt{2r}z}
+e^{-y-r+\sqrt{2r}z}-\mathcal{H}_{\mathbf{a},\alpha_{1},\alpha_{2}}^{\log(\mathcal{H}_{\alpha_{1}}\mathcal{H}_{\alpha_{2}})+x, \log
(\mathcal{H}_{a_{1},\alpha_{1}}\mathcal{H}_{a_{2},\alpha_{2}})+y}e^{-r+\sqrt{2r}z}\big)\right)\phi(z)dz
\end{eqnarray}
as $\mathbf{T}\rightarrow\infty,$ where
$$b_{\mathbf{a},\mathbf{T}}=a_{\mathbf{T}}+a_{\mathbf{T}}^{-1}\log\left((2\pi)^{-1/2}\mathcal{H}_{a_{1},\alpha_{1}}\mathcal{H}_{a_{2},\alpha_{2}}
(a_{\mathbf{T}})^{1/\alpha_{1}+1/\alpha_{2}-1/2}\right).$$ }

\textbf{Theorem 2.3}. {\sl Let $\{X(\mathbf{t}):\mathbf{t}\geq
\mathbf{0}\}$ be a homogeneous  Gaussian field with covariance
function $r(\mathbf{t})$ satisfying \textbf{A1}, \textbf{A2} and
\textbf{A3}. Then for any dense grids $\mathfrak{R}(p_{i})$,
$i=1,2$,
\begin{eqnarray}
\label{eq2.3}
&&P\left\{a_{\mathbf{T}}\big(M_{\mathbf{T}}-b_{\mathbf{T}}\big)\leq x, a_{\mathbf{T}}\big(M_{\mathbf{T}}^{\mathbf{p}}- b_{\mathbf{T}}\big)\leq y\right\}
\longrightarrow \int_{-\infty}^{+\infty}\exp\left(-e^{-\min(x,y)-r+\sqrt{2r}z}\right)\phi(z)dz
\end{eqnarray}
as $\mathbf{T}\rightarrow\infty.$ }

\textbf{Remark 2.2}.
i). In the literature, the Gaussian field $X(\mathbf{t})$ with correlation satisfying $\lim_{\mathbf{t}\rightarrow\infty}r(\mathbf{t})\log(t_{1}t_{2})=r\in[0,\infty)$ is called weakly and strongly dependent for $r=0$ and $r>0$, respectively,
see eg., Mittal and Ylvisaker (1975).
Theorems 2.1-2.3 show that for the weakly dependent case the two maxima are
asymptotically independent, dependent and coincide when the grid of the discrete time
points is a sparse grid, Pickands grid and dense grid, respectively. For the strongly dependent case,
the asymptotic independence between the two maxima does not hold anymore because of the
strong dependence of $X(\mathbf{t})$. However, in this case these two maxima are
asymptotically totally dependent if the grid of the discrete time
points is sufficiently dense, and asymptotically dependent if the
the grid points are sparse or Pickands grids. \\
ii). If $T_{1}=O(T_{2})$, as $\mathbf{T}\rightarrow\infty$, then
the condition that $r(0,t_{2})\log t_{2}$ and $r(t_{1},0)\log t_{1}$
are bounded in Assumption ${\bf A3}$ can be omitted.
Noting that Assumption ${\bf A3}$ is only used in the proof of Lemma 3.3,
it is easy to check this point from the bounds of $S_{\mathbf{T},22}$ and $M_{\mathbf{T},22}$ in the proofs of Lemma B1 and B3,
respectively. However, if $T_{1}=o(T_{2})$, as $\mathbf{T}\rightarrow\infty$,
then  Assumption ${\bf A3}$ can be weakened as: $\lim_{\mathbf{t}\rightarrow\infty}r(\mathbf{t})\log(t_{1}t_{2})=r\in[0,\infty)$ and $r(0,t_{2})\log t_{2}$ is bounded. A similar statement holds also for the case  $T_{2}=o(T_{1})$.

\section{Proofs}

First, define $\rho(\mathbf{T})=r/\log(T_{1}T_{2})$ and let $a> b$
be constants which will be determined in the proof of Lemma 3.3.
Following Piterbarg (2004), divide $[0,T_{i}]$ into intervals with
length $T_{i}^{a}$ alternating with shorter intervals with length
$T_{i}^{b}$, $i=1,2$. Note that the numbers of the long intervals is at most
$n_{i}=\lfloor T_{i}/(T_{i}^{a}+T_{i}^{b})\rfloor$, where $\lfloor x\rfloor$ denote the integral parts of $x$.
Let
$\mathbf{O_{i}}=
[(i_{1}-1)(T_{1}^{a}+T_{1}^{b}),(i_{1}-1)(T_{1}^{a}+T_{1}^{b})+T_{1}^{a}]\times[(i_{2}-1)(T_{2}^{a}+T_{2}^{b}),(i_{2}-1)(T_{2}^{a}+T_{2}^{b})+T_{2}^{a}]$,
$\mathbf{i}=\mathbf{1},\cdots,\mathbf{n}$ and
$\mathbf{O}=\cup_{\mathbf{i}} \mathbf{O}_{\mathbf{i}}$.
We will show blow that the remaining area $\mathbf{I_{T}}\backslash \mathbf{O}$  plays no role in our consideration.

Let $\{X_{\mathbf{i}}(\mathbf{t}), \mathbf{t}\geq \mathbf{0}\}$,
$\mathbf{i}\geq \mathbf{1}$ be independent copies of
$\{X(\mathbf{t}), \mathbf{t}\geq\mathbf{0}\}$ and
$\{\eta(\mathbf{t}), \mathbf{t}\geq \mathbf{0}\}$ be such that
$\eta(\mathbf{t})=X_{\mathbf{i}}(\mathbf{t})$ for $\mathbf{t}\in
\mathbf{E_{i}}:=
[(i_{1}-1)(T_{1}^{a}+T_{1}^{b}),i_{1}(T_{1}^{a}+T_{1}^{b}))\times[(i_{2}-1)(T_{2}^{a}+T_{2}^{b}),i_{2}(T_{2}^{a}+T_{2}^{b}))$, $\mathbf{i}=\mathbf{1},\cdots,\mathbf{n}$.
%Denote by $\gamma(\mathbf{s},\mathbf{t})$ the covariance function of
%$\{\eta(\mathbf{t}), \mathbf{t}\geq \mathbf{0}\}$. It is easy to see
%that
%\[
%  \gamma(\mathbf{s},\mathbf{t})=\left\{
% \begin{array}{cc}
%  {r(\mathbf{t},\mathbf{s})},    & \mathbf{s}\in \mathbf{E_{i}}, \mathbf{t}\in \mathbf{E_{j}}, \mathbf{i}=\mathbf{j};\\
%  {0},    & \mathbf{s}\in \mathbf{E_{i}}, \mathbf{t}\in \mathbf{E_{j}}, \mathbf{i}\neq\mathbf{j}.
% \end{array}
%  \right.
%\]\\
Define
$$\xi_{\mathbf{T}}(\mathbf{t})=\big(1-\rho(\mathbf{T})\big)^{1/2}\eta(\mathbf{t})+\rho^{1/2}(\mathbf{T})U, \ \ \mathbf{t}\in \mathbf{I_{T}},$$
where $U$ is a standard normal variable independent of
$\{\eta(\mathbf{t}), \mathbf{t}\geq \mathbf{0}\}$.  Denote by
$\varrho(\mathbf{s},\mathbf{t})$ the covariance function of
$\{\xi_{\mathbf{T}}(\mathbf{t}), \mathbf{t}\in \mathbf{I_{T}}\}$. It
is easy to check that
\[
  \varrho(\mathbf{s},\mathbf{t})=\left\{
 \begin{array}{cc}
  {r(\mathbf{t},\mathbf{s})+(1-r(\mathbf{t},\mathbf{s}))\rho(\mathbf{T})},    &\mathbf{s}\in \mathbf{E_{i}}, \mathbf{t}\in \mathbf{E_{j}}, \mathbf{i}=\mathbf{j};\\
  {\rho(\mathbf{T})},    & \mathbf{s}\in \mathbf{E_{i}}, \mathbf{t}\in \mathbf{E_{j}}, \mathbf{i}\neq\mathbf{j}.
 \end{array}
  \right.
\]\\
The proofs of our main results rely on the following Lemmas.
In the sequel, $C$ shall denote
positive constant whose values may vary from place to place.
% In the
%following part, let $C$ denote positive constants whose values may
%vary from place to place.

\textbf{Lemma 3.1}. {\sl Suppose that the grids
$\mathfrak{R}(p_{i})$, $i=1,2$ are  sparse grids or Pickands grids. For any
$B>0$, we have for all $x,y\in[-B,B]$,
\begin{eqnarray*}
&&\bigg|P\left\{a_{\mathbf{T}}\big(M_{\mathbf{T}}-b_{\mathbf{T}}\big)\leq x, a_{\mathbf{T}}\big(M_{\mathbf{T}}^{\mathbf{p}}- b_{\mathbf{T}}^{'}\big)\leq y\right\}\\
&&\ \ \ \ \ -P\left\{a_{\mathbf{T}}\big(\max_{\mathbf{t}\in \mathbf{O}}X(\mathbf{t})-b_{\mathbf{T}}\big)\leq x,
a_{\mathbf{T}}\big(\max_{\mathbf{t}\in\mathfrak{R}(\delta_{1})\times\mathfrak{R}(\delta_{2})\cap\mathbf{O}}X(\mathbf{t})- b_{\mathbf{T}}^{'}\big)\leq y\right\}\bigg|\rightarrow0
\end{eqnarray*}
as $\mathbf{T}\rightarrow\infty$, where
$b_{\mathbf{T}}^{'}=b_{\mathbf{T}}^{\mathbf{p}}$ for  sparse grids and
$b_{\mathbf{T}}^{'}=b_{\mathbf{a},\mathbf{T}}$ for  Pickands grids.
}

\textbf{Proof:} The proof is similar to that of Lemma 6 of Piterbarg
(2004). Clearly, we have
\begin{eqnarray}
\label{eq400}
&&\bigg|P\left\{a_{\mathbf{T}}\big(M_{\mathbf{T}}-b_{\mathbf{T}}\big)\leq x, a_{\mathbf{T}}\big(M_{\mathbf{T}}^{\mathbf{p}}- b_{\mathbf{T}}^{'}\big)\leq y\right\}\nonumber\\
&&\ \ \ \ \ -P\left\{a_{\mathbf{T}}\big(\max_{\mathbf{t}\in \mathbf{O}}X(\mathbf{t})-b_{\mathbf{T}}\big)\leq x,
a_{\mathbf{T}}\big(\max_{\mathbf{t}\in\mathfrak{R}(\delta_{1})\times\mathfrak{R}(\delta_{2})\cap\mathbf{O}}X(\mathbf{t})- b_{\mathbf{T}}^{'}\big)\leq y\right\}\bigg|\nonumber\\
&&\ \ \ \ \ \leq P\left\{\max_{\mathbf{t}\in \mathbf{I_{T}}\backslash\mathbf{O}}X(\mathbf{t})>b_{\mathbf{T}}+x/a_{\mathbf{T}}\right\}
+P\left\{\max_{\mathbf{t}\in\mathfrak{R}(\delta_{1})\times\mathfrak{R}(\delta_{2})\cap\mathbf{I_{T}}\backslash\mathbf{O}}X(\mathbf{t})>b_{\mathbf{T}}^{'}+y/a_{\mathbf{T}}\right\}
\end{eqnarray}
By Theorem 7.2 of Piterbarg (1996) (denote by $mes(\cdot)$ the Lebesgue measure)
\begin{eqnarray*}
P\left\{\max_{\mathbf{t}\in \mathbf{I_{T}}\backslash\mathbf{O}}X(\mathbf{t})>b_{\mathbf{T}}+x/a_{\mathbf{T}}\right\}
&=&O(1)mes(\mathbf{I_{T}}\backslash\mathbf{O})(b_{\mathbf{T}}+x/a_{\mathbf{T}})^{2/\alpha_{1}+2/\alpha_{2}}\Psi(b_{\mathbf{T}}+x/a_{\mathbf{T}})\\
&=&O(1)\frac{mes(\mathbf{I_{T}}\backslash\mathbf{O})}{T_{1}T_{2}}\\
&\leq &O(1)\frac{n_{1}T_{1}^{b}(T_{2}^{a}+T_{2}^{b})+n_{2}T_{2}^{b}(T_{1}^{a}+T_{1}^{b})}{T_{1}T_{2}}\rightarrow0
\end{eqnarray*}
as $\mathbf{T}\rightarrow\infty$, by the choice of $a_{\mathbf{T}}$
and $b_{\mathbf{T}}$. In light of (\ref{eq313}) and
(\ref{eq319}) in the Appendix for a sparse grid and Pickands grid, respectively, we can
get the same estimation for the second probability in the right-hand side
of (\ref{eq400}), hence the proof is complete. \hfill $\Box$

For the proofs we need also the following auxiliary grids $\mathfrak{R}(q_{i})$
with $q_{i}=\gamma_{i}(2\log T_{1}T_{2})^{-1/\alpha_{i}}$ and $\gamma_{i}>0$, $i=1,2$.

\textbf{Lemma 3.2}. {\sl Suppose that the grids
$\mathfrak{R}(p_{i})$, $i=1,2$ are sparse grids or Pickands grids. For any
$B>0$, we have for all $x,y\in[-B,B]$
\begin{eqnarray*}
&&\bigg|P\left\{a_{\mathbf{T}}\big(\max_{\mathbf{t}\in \mathbf{O}}X(\mathbf{t})-b_{\mathbf{T}}\big)\leq x,
a_{\mathbf{T}}\big(\max_{\mathbf{t}\in\mathfrak{R}(p_{1})\times\mathfrak{R}(p_{2})\cap\mathbf{O}}X(\mathbf{t})- b_{\mathbf{T}}^{'}\big)\leq y\right\}\\
&&\ \ \ \ \ -P\left\{a_{\mathbf{T}}\big(\max_{\mathbf{t}\in \mathfrak{R}(q_{1})\times\mathfrak{R}(q_{2})\cap\mathbf{O}}X(\mathbf{t})-b_{T}\big)\leq x,
     a_{\mathbf{T}}\big(\max_{\mathbf{t}\in\mathfrak{R}(p_{1})\times\mathfrak{R}(p_{2})\cap\mathbf{O}}X(\mathbf{t})- b_{\mathbf{T}}^{'}\big)\leq y\right\}\bigg|\rightarrow0
\end{eqnarray*}
as $\mathbf{T}\rightarrow\infty$ and $\gamma_{i}\downarrow0$, where
$b_{\mathbf{T}}^{'}=b_{\mathbf{T}}^{\mathbf{p}}$ for  sparse grids and
$b_{\mathbf{T}}^{'}=b_{\mathbf{a},\mathbf{T}}$ for  Pickands grids.
}

\textbf{Proof:} It follows from Lemma A4.\hfill$\Box$

The following lemma plays a crucial role in the proofs of Theorems 2.1 and 2.2.

\textbf{Lemma 3.3}. {\sl Suppose that the grids
$\mathfrak{R}(p_{i})$, $i=1,2$ are  sparse grids or  Pickands grids. For any
$B>0$ we have for all $x,y\in[-B,B]$,
\begin{eqnarray*}
&&\bigg|P\left\{a_{\mathbf{T}}\big(\max_{\mathbf{t}\in \mathfrak{R}(q_{1})\times\mathfrak{R}(q_{2})\cap\mathbf{O}}X(\mathbf{t})-b_{T}\big)\leq x,
     a_{\mathbf{T}}\big(\max_{\mathbf{t}\in\mathfrak{R}(p_{1})\times\mathfrak{R}(p_{2})\cap\mathbf{O}}X(\mathbf{t})- b_{\mathbf{T}}^{'}\big)\leq y\right\}\\
&&\ \ \ \ \ -P\left\{a_{\mathbf{T}}\big(\max_{\mathbf{t}\in \mathfrak{R}(q_{1})\times\mathfrak{R}(q_{2})\cap\mathbf{O}}\xi_{\mathbf{T}}(\mathbf{t})-b_{\mathbf{T}}\big)\leq x,
     a_{\mathbf{T}}\big(\max_{\mathbf{t}\in\mathfrak{R}(p_{1})\times\mathfrak{R}(p_{2})\cap\mathbf{O}}\xi_{\mathbf{T}}(\mathbf{t})- b_{\mathbf{T}}^{'}\big)\leq y\right\}\bigg|\rightarrow 0
\end{eqnarray*}
as $\mathbf{T}\rightarrow\infty$, where
$b_{\mathbf{T}}^{'}=b_{\mathbf{T}}^{\mathbf{p}}$ for  sparse grids and
$b_{\mathbf{T}}^{'}=b_{\mathbf{a},\mathbf{T}}$ for  Pickands grids.
}

\textbf{Proof:} For the sake of simplicity, let $u_{\mathbf{T}}=
b_{\mathbf{T}}+x/a_{\mathbf{T}}$, $u_{\mathbf{T}}'=
b_{\mathbf{T}}'+y/a_{\mathbf{T}}$. Using the Normal
Comparison Lemma (see eg. Leadbetter et al. (1983) and Piterbarg (1996)),
we have
\begin{eqnarray*}
\label{eq401}
&&\bigg|P\left\{a_{\mathbf{T}}\big(\max_{\mathbf{t}\in \mathfrak{R}(q_{1})\times\mathfrak{R}(q_{2})\cap\mathbf{O}}X(\mathbf{t})-b_{\mathbf{T}}\big)\leq x,
     a_{\mathbf{T}}\big(\max_{\mathbf{t}\in\mathfrak{R}(p_{1})\times\mathfrak{R}(p_{2})\cap\mathbf{O}}X(\mathbf{t})- b_{\mathbf{T}}^{'}\big)\leq y\right\}\nonumber\\
&&\ \ \ \ \ -P\left\{a_{\mathbf{T}}\big(\max_{\mathbf{t}\in \mathfrak{R}(q_{1})\times\mathfrak{R}(q_{2})\cap\mathbf{O}}\xi_{\mathbf{T}}(\mathbf{t})-b_{\mathbf{T}}\big)\leq x,
     a_{\mathbf{T}}\big(\max_{\mathbf{t}\in\mathfrak{R}(p_{1})\times\mathfrak{R}(p_{2})\cap\mathbf{O}}\xi_{\mathbf{T}}(\mathbf{t})- b_{\mathbf{T}}^{'}\big)\leq y\right\}\bigg|\nonumber\\
&&\leq \sum_{\mathbf{kq}\in \mathbf{O_{i}},\mathbf{lq}\in \mathbf{O_{j}}\atop \mathbf{kq}\neq \mathbf{lq}, \mathbf{1}\leq \mathbf{i,j}\leq \mathbf{n} }|r(\mathbf{kq},\mathbf{lq})-\varrho(\mathbf{kq},\mathbf{lq})|
\int_{0}^{1}\frac{1}{\sqrt{1-r^{(h)}(\mathbf{kq},\mathbf{lq})}}\exp\left(-\frac{u_{\mathbf{T}}^{2}}{1+r^{(h)}(\mathbf{kq},\mathbf{lq})}\right)dh\nonumber\\
&&+\sum_{\mathbf{\mathbf{kp}}\in \mathbf{O_{i}},\mathbf{\mathbf{lp}}\in \mathbf{O_{j}}\atop \mathbf{kp}\neq \mathbf{lp}, \mathbf{1}\leq \mathbf{i,j}\leq \mathbf{n} }|r(\mathbf{kp},\mathbf{lp})-\varrho(\mathbf{kp},\mathbf{lp})|
\int_{0}^{1}\frac{1}{\sqrt{1-r^{(h)}(\mathbf{kp},\mathbf{lp})}}\exp\left(-\frac{u_{\mathbf{T}}'^{2}}{1+r^{(h)}(\mathbf{kp},\mathbf{lp})}\right)dh\nonumber\\
&&+\sum_{\mathbf{kq}\in \mathbf{O_{i}},\mathbf{lp}\in \mathbf{O_{j}}\atop \mathbf{kq}\neq \mathbf{lp}, \mathbf{1}\leq \mathbf{i,j}\leq \mathbf{n} }|r(\mathbf{kq},\mathbf{lp})-\varrho(\mathbf{kq},\mathbf{lp})|
\int_{0}^{1}\frac{1}{\sqrt{1-r^{(h)}(\mathbf{kq},\mathbf{lp})}}\exp\left(-\frac{u_{\mathbf{T}}'^{2}+u_{\mathbf{T}}^{2}}{2(1+r^{(h)}(\mathbf{kq},\mathbf{lp}))}\right)dh,
\end{eqnarray*}
where
$r^{(h)}(\mathbf{kq},\mathbf{lq})=hr(\mathbf{kq},\mathbf{lq})+(1-h)\varrho(\mathbf{kq},\mathbf{lq})$.
Now, the lemma follows from Lemmas B1-B3 in the Appendix B.
\hfill$\Box$

\textbf{Lemma 3.4}. {\sl Suppose that the grids
$\mathfrak{R}(p_{i})$, $i=1,2$ are sparse grids or  Pickands grids. For any
$B>0$ we have for all $x,y\in[-B,B]$ and the grids $\mathfrak{R}(q_{i})$
with $q_{i}=\gamma_{i}(2\log T_{1}T_{2})^{-1/\alpha_{i}}$ and $\gamma_{i}>0$, $i=1,2$
\begin{eqnarray*}
&&\bigg|P\left\{a_{\mathbf{T}}\big(\max_{\mathbf{t}\in \mathfrak{R}(q_{1})\times\mathfrak{R}(q_{2})\cap\mathbf{O}}\xi_{\mathbf{T}}(\mathbf{t})-b_{\mathbf{T}}\big)\leq x,
     a_{\mathbf{T}}\big(\max_{\mathbf{t}\in\mathfrak{R}(p_{1})\times\mathfrak{R}(p_{2})\cap\mathbf{O}}\xi_{\mathbf{T}}(\mathbf{t})- b_{\mathbf{T}}^{'}\big)\leq y\right\}\\
&&\ \ \ \ \ \ -\int_{-\infty}^{+\infty}\prod_{i_{1}=1}^{n_{1}}\prod_{i_{2}=1}^{n_{2}}
P\left\{\max_{\mathbf{t}\in \mathbf{O_{i}}}\eta(\mathbf{t})\leq u_{\mathbf{T}}^{*},
   \max_{\mathbf{t}\in\mathbf{O_{i}}}\eta(\mathbf{t})\leq u_{\mathbf{T}}^{*'}\right\}\phi(z)dz\bigg|\rightarrow 0,
\end{eqnarray*}
as $\gamma_{i}\downarrow 0$, where
\begin{eqnarray}
\label{equ413}
u_{\mathbf{T}}^{*}:=\frac{b_{\mathbf{T}}+x/a_{\mathbf{T}}-\rho^{1/2}(\mathbf{T})z}{(1-\rho(\mathbf{T}))^{1/2}}
=\frac{x+r-\sqrt{2r}z}{a_{\mathbf{T}}}+b_{\mathbf{T}}+o(a_{\mathbf{T}}^{-1}),
\end{eqnarray}
and
\begin{eqnarray}
\label{equ414}
u_{\mathbf{T}}^{*'}:=\frac{b_{\mathbf{T}}^{'}+y/a_{\mathbf{T}}-\rho^{1/2}(\mathbf{T})z}{(1-\rho(\mathbf{T}))^{1/2}}
=\frac{y+r-\sqrt{2r}z}{a_{\mathbf{T}}}+b_{\mathbf{T}}^{'}+o(a_{\mathbf{T}}^{-1})
\end{eqnarray}
with $b_{\mathbf{T}}^{'}=b_{\mathbf{T}}^{\mathbf{p}}$ for  sparse grids and
$b_{\mathbf{T}}^{'}=b_{\mathbf{a},\mathbf{T}}$ for Pickands grids.
}

\textbf{Proof:} First, by the definition of
$\{\xi_{\mathbf{T}}(\mathbf{t}), \mathbf{0}\leq \mathbf{t}\leq
\mathbf{T}\}$, we have
\begin{eqnarray}
\label{eq412}
&&P\left\{a_{\mathbf{T}}\big(\max_{\mathbf{t}\in \mathfrak{R}(q_{1})\times\mathfrak{R}(q_{2})\cap\mathbf{O}}\xi_{\mathbf{T}}(\mathbf{t})-b_{\mathbf{T}}\big)\leq x,
     a_{\mathbf{T}}\big(\max_{\mathbf{t}\in\mathfrak{R}(p_{1})\times\mathfrak{R}(p_{2})\cap\mathbf{O}}\xi_{\mathbf{T}}(\mathbf{t})- b_{\mathbf{T}}^{'}\big)\leq y\right\}\nonumber\\
&&\ \ =\int_{-\infty}^{+\infty}P\left\{\max_{\mathbf{t}\in \mathfrak{R}(q_{1})\times\mathfrak{R}(q_{2})\cap\mathbf{O}}\eta(\mathbf{t})\leq u_{\mathbf{T}}^{*},
     \max_{\mathbf{t}\in\mathfrak{R}(p_{1})\times\mathfrak{R}(p_{2})\cap\mathbf{O}}\eta(\mathbf{t})\leq u_{\mathbf{T}}^{*'}\right\}\phi(z)dz\nonumber\\
&&\ \ =\int_{-\infty}^{+\infty}\prod_{i_{1}=1}^{n_{1}}\prod_{i_{2}=1}^{n_{2}}
P\left\{\max_{\mathbf{t}\in \mathfrak{R}(q_{1})\times\mathfrak{R}(q_{2})\cap \mathbf{O_{i}}}\eta(\mathbf{t})\leq u_{\mathbf{T}}^{*},
   \max_{t\in\mathfrak{R}(p_{1})\times\mathfrak{R}(p_{2})\cap\mathbf{O_{i}}}\eta(\mathbf{t})\leq u_{\mathbf{T}}^{*'}\right\}\phi(z)dz.
\end{eqnarray}
As for the discrete case, see page 137 on Leadbetter et al. (1983), a direct calculation leads to
\begin{eqnarray*}
u_{\mathbf{T}}^{*}&=&\frac{x+r-\sqrt{2r}z}{a_{\mathbf{T}}}+b_{\mathbf{T}}+o(a_{\mathbf{T}}^{-1})
\end{eqnarray*}
and
\begin{eqnarray*}
u_{\mathbf{T}}^{*'}&=&\frac{y+r-\sqrt{2r}z}{a_{\mathbf{T}}}+b_{\mathbf{T}}^{'}+o(a_{\mathbf{T}}^{-1}).
\end{eqnarray*}
Next, by Lemma A4  and the dominated convergence theorem, we have
\begin{eqnarray}
\label{eq413}
&&\bigg|\int_{-\infty}^{+\infty}\prod_{i_{1}=1}^{n_{1}}\prod_{i_{2}=1}^{n_{2}}
P\left\{\max_{\mathbf{t}\in \mathfrak{R}(q_{1})\times\mathfrak{R}(q_{2})\cap \mathbf{O_{i}}}\eta(\mathbf{t})\leq u_{\mathbf{T}}^{*},
   \max_{\mathbf{t}\in\mathfrak{R}(p_{1})\times\mathfrak{R}(p_{2})\cap\mathbf{O_{i}}}\eta(\mathbf{t})\leq u_{\mathbf{T}}^{*'}\right\}\phi(z)dz\nonumber\\
&&\ \ \ \ \ \ -\int_{-\infty}^{+\infty}\prod_{i_{1}=1}^{n_{1}}\prod_{i_{2}=1}^{n_{2}}
P\left\{\max_{\mathbf{t}\in \mathbf{O_{i}}}\eta(\mathbf{t})\leq u_{\mathbf{T}}^{*},
   \max_{\mathbf{t}\in\mathfrak{R}(p_{1})\times\mathfrak{R}(p_{2})\cap\mathbf{O_{i}}}\eta(\mathbf{t})\leq u_{\mathbf{T}}^{*'}\right\}\phi(z)dz\bigg|\rightarrow 0
\end{eqnarray}
as $\gamma_{i}\downarrow 0$. Lemma 3.4 now follows by combining
(\ref{eq412}) with (\ref{eq413}). \hfill$\Box$

\textbf{Proof of Theorem 2.1}. From Lemmas 3.1-3.4, we known that in order to prove Theorem 2.1, it
suffice to show that
\begin{eqnarray*}
&&\bigg|\int_{-\infty}^{+\infty}\prod_{i_{1}=1}^{n_{1}}\prod_{i_{2}=1}^{n_{2}}
P\left\{\max_{\mathbf{t}\in \mathbf{O_{i}}}\eta(\mathbf{t})\leq u_{\mathbf{T}}^{*},
   \max_{\mathbf{t}\in\mathfrak{R}(p_{1})\times\mathfrak{R}(p_{2})\cap\mathbf{O_{i}}}\eta(\mathbf{t})\leq u_{\mathbf{T}}^{*'}\right\}\phi(z)dz\\
&&\ \ \ \ \ \ \ -\int_{-\infty}^{+\infty}\exp\left(-\big(e^{-x-r+\sqrt{2r}z}+e^{-y-r+\sqrt{2r}z}\big)\right)\phi(z)dz\bigg|\rightarrow 0
\end{eqnarray*}
as $\mathbf{T}\rightarrow\infty$, where $u_{\mathbf{T}}^{*}$ and
$u_{\mathbf{T}}^{*'}$ are defined in Lemma 3.4.
Using the homogeneity of $\{\eta(\mathbf{t}), \mathbf{t}\geq \mathbf{0}\}$,
\begin{eqnarray*}
&&\prod_{i_{1}=1}^{n_{1}}\prod_{i_{2}=1}^{n_{2}}
P\left\{\max_{\mathbf{t}\in \mathbf{O_{i}}}\eta(\mathbf{t})\leq u_{\mathbf{T}}^{*},
   \max_{\mathbf{t}\in\mathfrak{R}(p_{1})\times\mathfrak{R}(p_{2})\cap\mathbf{O_{i}}}\eta(\mathbf{t})\leq u_{\mathbf{T}}^{*'}\right\}\\
&&=\left(P\left\{\max_{\mathbf{t}\in [0,T_{1}^{a}]\times[0,T_{2}^{a}]}\eta(\mathbf{t})\leq u_{\mathbf{T}}^{*},
   \max_{\mathbf{t}\in\mathfrak{R}(p_{1})\times\mathfrak{R}(p_{2})\cap [0,T_{1}^{a}]\times[0,T_{2}^{a}]}\eta(\mathbf{t})\leq u_{\mathbf{T}}^{*'}\right\}\right)^{n_{1}n_{2}}\\
&&=\exp\left(n_{1}n_{2}\log\left(P\left\{\max_{\mathbf{t}\in [0,T_{1}^{a}]\times[0,T_{2}^{a}]}\eta(\mathbf{t})\leq u_{\mathbf{T}}^{*},
   \max_{\mathbf{t}\in\mathfrak{R}(p_{1})\times\mathfrak{R}(p_{2})\cap [0,T_{1}^{a}]\times[0,T_{2}^{a}]}\eta(\mathbf{t})\leq u_{\mathbf{T}}^{*'}\right\}\right)\right)\\
&&=\exp\left(-n_{1}n_{2}\left(1-P\left\{\max_{\mathbf{t}\in [0,T_{1}^{a}]\times[0,T_{2}^{a}]}\eta(\mathbf{t})\leq u_{\mathbf{T}}^{*},
   \max_{\mathbf{t}\in\mathfrak{R}(p_{1})\times\mathfrak{R}(p_{2})\cap [0,T_{1}^{a}]\times[0,T_{2}^{a}]}\eta(\mathbf{t})\leq u_{\mathbf{T}}^{*'}\right\}\right)+R_{\mathbf{n}}\right).
\end{eqnarray*}
Since
$$P_{\mathbf{n}}=:P\left\{\max_{\mathbf{t}\in [0,T_{1}^{a}]\times[0,T_{2}^{a}]}\eta(\mathbf{t})\leq u_{\mathbf{T}}^{*},
   \max_{\mathbf{t}\in\mathfrak{R}(p_{1})\times\mathfrak{R}(p_{2})\cap [0,T_{1}^{a}]\times[0,T_{2}^{a}]}\eta(\mathbf{t})\leq u_{\mathbf{T}}^{*'}\right\}\rightarrow 1,$$
as $\mathbf{T}\rightarrow\infty$, we get that the remainder
$R_{\mathbf{n}}$ can be estimated as
$R_{\mathbf{n}}=o(n_{1}n_{2}(1-P_{\mathbf{n}}))$. Using Lemma A2,
 (\ref{equ413}) and (\ref{equ414}), we get that
\begin{eqnarray*}
&&n_{1}n_{2}\left(1-P\left\{\max_{\mathbf{t}\in [0,T_{1}^{a}]\times[0,T_{2}^{a}]}\eta(\mathbf{t})\leq u_{\mathbf{T}}^{*},
   \max_{\mathbf{t}\in\mathfrak{R}(p_{1})\times\mathfrak{R}(p_{2})\cap [0,T_{1}^{a}]\times[0,T_{2}^{a}]}\eta(\mathbf{t})\leq u_{\mathbf{T}}^{*'}\right\}\right)\\
&&\thicksim n_{1}n_{2}T_{1}^{a}T_{2}^{a}T_{1}^{-1}T_{2}^{-1}\left(e^{-x-r+\sqrt{2r}z}+e^{-y-r+\sqrt{2r}z}\right)\\
&&\thicksim e^{-x-r+\sqrt{2r}z}+e^{-y-r+\sqrt{2r}z},
\end{eqnarray*}
which combined with the dominated convergence theorem completes the
proof of Theorem 2.1.\hfill$\Box$

\textbf{Proof of Theorem 2.2}. The proof of the first assertion of Theorem 2.2 can be found in Subsection 4.2.
Next, we give the proof of the second assertion. In view of Lemmas 3.1-3.4 in order to establish the proof we
need to show
\begin{eqnarray*}
&&\bigg|\int_{-\infty}^{+\infty}\prod_{i_{1}=1}^{n_{1}}\prod_{i_{2}=1}^{n_{2}}
P\left\{\max_{\mathbf{t}\in \mathbf{O_{i}}}\eta(\mathbf{t})\leq u_{\mathbf{T}}^{*},
   \max_{t\in\mathfrak{R}(p_{1})\times\mathfrak{R}(p_{2})\cap\mathbf{O_{i}}}\eta(\mathbf{t})\leq u_{\mathbf{T}}^{*'}\right\}\phi(z)dz\\
&&\ \ \ \ \ \ \ -\int_{-\infty}^{+\infty}\exp\left(-\big(e^{-x-r+\sqrt{2r}z}+e^{-y-r+\sqrt{2r}z}
       -H_{a,\alpha}^{\log H_{a,\alpha}+x,\log H_{\alpha}+y}e^{-r+\sqrt{2r}z}\big)\right)\phi(z)dz\bigg|\rightarrow 0
\end{eqnarray*}
as $\mathbf{T}\rightarrow\infty$, where $u_{\mathbf{T}}^{*}$ and
$u_{\mathbf{T}}^{*'}$ are defined in Lemma 3.4. Similar to the proof
of Theorem 2.1, using Lemma A3, we get
\begin{eqnarray*}
&&n_{1}n_{2}\left(1-P\left\{\max_{\mathbf{t}\in [0,T_{1}^{a}]\times[0,T_{2}^{a}]}\eta(\mathbf{t})\leq u_{\mathbf{T}}^{*},
   \max_{\mathbf{t}\in\mathfrak{R}(p_{1})\times\mathfrak{R}(p_{2})\cap [0,T_{1}^{a}]\times[0,T_{2}^{a}]}\eta(\mathbf{t})\leq u_{\mathbf{T}}^{*'}\right\}\right)\\
&&=n_{1}n_{2}P\left\{\max_{\mathbf{t}\in [0,T_{1}^{a}]\times[0,T_{2}^{a}]}\eta(\mathbf{t})> u_{\mathbf{T}}^{*}\right\}
   +n_{1}n_{2}P\left\{\max_{\mathbf{t}\in\mathfrak{R}(p_{1})\times\mathfrak{R}(p_{2})\cap [0,T_{1}^{a}]\times[0,T_{2}^{a}]}\eta(\mathbf{t})> u_{\mathbf{T}}^{*'}\right\}\\
&&\ \ -n_{1}n_{2}P\left\{\max_{\mathbf{t}\in [0,T_{1}^{a}]\times[0,T_{2}^{a}]}\eta(\mathbf{t})> u_{\mathbf{T}}^{*},
   \max_{\mathbf{t}\in\mathfrak{R}(p_{1})\times\mathfrak{R}(p_{2})\cap [0,T_{1}^{a}]\times[0,T_{2}^{a}]}\eta(\mathbf{t})> u_{\mathbf{T}}^{*'}\right\}\\
&&\thicksim n_{1}n_{2}T_{1}^{a}T_{2}^{a}T_{1}^{-1}T_{2}^{-1}\left(e^{-x-r+\sqrt{2r}z}+e^{-y-r+\sqrt{2r}z}\right)\\
&&\ \ -n_{1}n_{2}P\left\{\max_{\mathbf{t}\in [0,T_{1}^{a}]\times[0,T_{2}^{a}]}\eta(\mathbf{t})> u_{\mathbf{T}}^{*},
   \max_{\mathbf{t}\in\mathfrak{R}(p_{1})\times\mathfrak{R}(p_{2})\cap [0,T_{1}^{a}]\times[0,T_{2}^{a}]}\eta(\mathbf{t})> u_{\mathbf{T}}^{*'}\right\},
\end{eqnarray*}
as $\mathbf{T}\rightarrow\infty$. To transform the last term, using
(\ref{equ413}) and (\ref{equ414}), we get
\begin{eqnarray*}
u_{\mathbf{T}}^{*}%&=&\frac{b_{\mathbf{T}}+x/a_{\mathbf{T}}-\rho^{1/2}(\mathbf{T})z}{(1-\rho(\mathbf{T}))^{1/2}}\\
&=&\frac{x+r-\sqrt{2r}z}{a_{\mathbf{T}}}+b_{\mathbf{T}}+o(a_{\mathbf{T}}^{-1})\\
&=&u_{\mathbf{T}}^{*'}+b_{\mathbf{T}}-b_{a,\mathbf{T}}+(x-y)/a_{\mathbf{T}}\\
&=&u_{\mathbf{T}}^{*'}+\frac{\log (\mathcal{H}_{\alpha_{1}}\mathcal{H}_{\alpha_{2}})-\log (\mathcal{H}_{a_{1},\alpha_{1}}\mathcal{H}_{a_{2},\alpha_{2}})+x-y}{u_{\mathbf{T}}^{*'}}+O\left((\log\log (T_{1}T_{2}))^{2}(\log T_{1}T_{2})^{-3/2}\right).
\end{eqnarray*}
Observing that $u_{\mathbf{T}}^{*'}\thicksim (2\log
T_{1}T_{2})^{1/2}$, we see that the reminder $O(\cdot)$ plays a
negligible role. Therefore, by (\ref{eq318}) in Appendix A
\begin{eqnarray*}
&&n_{1}n_{2}P\left\{\max_{\mathbf{t}\in [0,T_{1}^{a}]\times[0,T_{2}^{a}]}\eta(\mathbf{t})> u_{\mathbf{T}}^{*},
   \max_{\mathbf{t}\in\mathfrak{R}(p_{1})\times\mathfrak{R}(p_{2})\cap [0,T_{1}^{a}]\times[0,T_{2}^{a}]}\eta(\mathbf{t})> u_{\mathbf{T}}^{*'}\right\}\\
&&=n_{1}n_{2}T_{1}^{a}T_{2}^{a}\mathcal{H}_{\mathbf{a},\alpha_{1},\alpha_{2}}^{0,Z_{x,y}}(u_{T}^{*'})^{2/\alpha_{1}+2/\alpha_{2}}\Psi(u_{T}^{*'})(1+o(1))\\
&&=n_{1}n_{2}T_{1}^{a}T_{2}^{a}(T_{1}T_{2})^{-1}\mathcal{H}_{\mathbf{a},\alpha_{1},\alpha_{2}}^{0,Z_{x,y}}(\mathcal{H}_{a_{1},\alpha_{1}}\mathcal{H}_{a_{2},\alpha_{2}})^{-1}e^{-y-r+\sqrt{2r}z}(1+o(1)),
\end{eqnarray*}
where $Z_{x,y}=\log
(\mathcal{H}_{\alpha_{1}}\mathcal{H}_{\alpha_{2}})-\log
(\mathcal{H}_{a_{1},\alpha_{1}}\mathcal{H}_{a_{2},\alpha_{2}})+x-y$.
Next, changing the variables in the definition of
$\mathcal{H}_{\mathbf{a},\alpha_{1},\alpha_{2}}^{x,y}$ we get that
$\mathcal{H}_{\mathbf{a},\alpha_{1},\alpha_{2}}^{0,Z_{x,y}}(\mathcal{H}_{a_{1},\alpha_{1}}\mathcal{H}_{a_{2},\alpha_{2}})^{-1}e^{-y}
=\mathcal{H}_{\mathbf{a},\alpha_{1},\alpha_{2}}^{\log
(\mathcal{H}_{\alpha_{1}}\mathcal{H}_{\alpha_{2}})+x, \log
(\mathcal{H}_{a_{1},\alpha_{1}}\mathcal{H}_{a_{2},\alpha_{2}})+y}$.
This and the dominated convergence theorem conclude the proof of
Theorem 2.2.\hfill$\Box$

\textbf{Proof of Theorem 2.3}. In view of Lemma A4 we have
\begin{eqnarray*}
&&\bigg|P\left\{a_{\mathbf{T}}\big(M_{\mathbf{T}}-b_{\mathbf{T}}\big)\leq x, a_{\mathbf{T}}\big(M_{\mathbf{T}}^{\mathbf{p}}- b_{\mathbf{T}}\big)\leq y\right\}
-P\left\{a_{\mathbf{T}}\big(M_{\mathbf{T}}-b_{\mathbf{T}}\big)\leq x, a_{\mathbf{T}}\big(M_{\mathbf{T}}- b_{\mathbf{T}}\big)\leq y\right\}\bigg|\\
&&\leq\bigg|P\left\{a_{\mathbf{T}}\big(M_{\mathbf{T}}^{\mathbf{p}}- b_{\mathbf{T}}\big)\leq y\right\}
-P\left\{a_{\mathbf{T}}\big(M_{\mathbf{T}}- b_{\mathbf{T}}\big)\leq y\right\}\bigg| \rightarrow 0, \quad \mathbf{T}\rightarrow\infty.
\end{eqnarray*}
Next, applying  Corollary 2.1, we get
\begin{eqnarray*}
P\left\{a_{\mathbf{T}}\big(M_{\mathbf{T}}-b_{\mathbf{T}}\big)\leq x, a_{\mathbf{T}}\big(M_{\mathbf{T}}- b_{\mathbf{T}}\big)\leq y\right\}
&=&P\{a_{\mathbf{T}}\big(M_{\mathbf{T}}-b_{\mathbf{T}}\big)\leq \min(x, y)\}\\
&\rightarrow& \int_{-\infty}^{+\infty}\exp\left(-e^{-\min(x,y)-r+\sqrt{2r}z}\right)\phi(z)dz, \quad \mathbf{T}\rightarrow\infty,
\end{eqnarray*}
hence the proof is complete. \hfill$\Box$

\section{Appendix A}

In this section, we give some auxiliary results, which extend Lemmas 1-4 of Piterbarg (2004) from stationary Gaussian processes to  Gaussian random fields.
The ideas of the proofs are very close to that of the above mentioned lemmas. In the following subsections, we suppose that Assumptions \textbf{A1} and \textbf{A2} hold.

\subsection{Sparse grid}
In this subsection, suppose $\mathfrak{R}(p_{i})$, $i=1,2$ are sparse grids. We will
use the notations $u=\sqrt{2\log T_{1}T_{2}}$, so that
$p_{i}=p_{i}(u)=l_{i}(u)u^{-2/\alpha_{i}}$,
$i=1,2$, where $l_{i}(u)\rightarrow\infty$ as
$u\rightarrow\infty$, with $p_{i}(u)\leq p_{0}$
for some positive $p_{0}$, in particular,
$p_{i}(u)=p_{0}$. Let
$\mathbf{L_{p}}=[-p_{1},p_{1}]\times[-p_{2},p_{2}]$.
First  we consider the following probability
$$P(u,x)=P\left(X(\mathbf{0})>u, \max_{\mathbf{t}\in\mathbf{L_{p}}}X(\mathbf{t})>
u+\frac{(\frac{2}{\alpha_{1}}+\frac{2}{\alpha_{2}})\log
u+\log(p_{1}p_{2})+x}{u}\right),$$ where $x$ is varies in
a closed interval, say, $x\in [-A,A]$ with  $A<\infty$. For simplicity, we denote
$$v:=\sqrt{(\frac{2}{\alpha_{1}}+\frac{2}{\alpha_{2}})\log u+\log(p_{1}p_{2})}.$$
By (\ref{eq2.0}) (see also Theorem 7.1 of Piterbarg (1996)), we have
\begin{eqnarray}
\label{eq311}
P\left(\max_{\mathbf{t}\in\mathbf{L_{p}}}X(\mathbf{t})>u+\frac{v^{2}+x}{u}\right)
&=& 4p_{1}p_{2}\mathcal{H}_{\alpha_{1}}\mathcal{H}_{\alpha_{2}}
            \left(u+\frac{v^{2}+x}{u}\right)^{2/\alpha_{1}+2/\alpha_{2}}\Psi\left(u+\frac{v^{2}+x}{u}\right)(1+o(1))\nonumber\\
%&&\thicksim 4p_{1}p_{2}\mathcal{H}_{\alpha_{1}}\mathcal{H}_{\alpha_{2}}e^{-x}
%            \left(u\right)^{2/\alpha_{1}+2/\alpha_{2}}\exp\left(-(\frac{2}{\alpha_{1}}+\frac{2}{\alpha_{2}})\log u-\log(p_{1}p_{2})\right)\Psi(u)\nonumber\\
&=& 4\mathcal{H}_{\alpha_{1}}\mathcal{H}_{\alpha_{2}}e^{-x}\Psi(u)(1+o(1))
\end{eqnarray}
as $u\rightarrow\infty$.

\textbf{Lemma A1}. {\sl We have $P(u,x)=o(\Psi(u))$ as
$u\rightarrow\infty$.}

\textbf{Proof:} Write
$w=\left((\frac{2}{\alpha_{1}}+\frac{2}{\alpha_{2}})\log
u+\log(p_{1}p_{2})+x\right)/u$. We have,
$$P(u,x)\leq P\left(\max_{\mathbf{t}\in\mathbf{L_{p}}}(X(\mathbf{0})+X(\mathbf{t}))>2u+w\right).$$
Let $p_{1}', p_{2}'$ be so small that $1-r(\mathbf{t})\leq
2(|t_{1}|^{\alpha_{1}}+|t_{2}|^{\alpha_{2}})$ for all $\mathbf{t}\in
\mathbf{I_{p'}}= [-p_{1}',
p_{1}']\times[-p_{2}', p_{2}']$. If
$\mathbf{L_{p}}\cap\mathbf{L_{p'}}\neq \emptyset$, write
\begin{eqnarray}
\label{eq312}
P\left(\max_{\mathbf{t}\in\mathbf{L_{p}}}(X(\mathbf{0})+X(\mathbf{t}))>2u+w\right)
&\leq& P\left(\max_{\mathbf{t}\in\mathbf{L_{p'}}}(X(\mathbf{0})+X(\mathbf{t}))>2u+w\right)\nonumber\\
&+&P\left(\max_{\mathbf{t}\in \mathbf{L_{p}}\diagdown \mathbf{L_{p'}}}(X(\mathbf{0})+X(\mathbf{t}))>2u+w\right).
\end{eqnarray}
The variance of the field $X(\mathbf{0})+X(\mathbf{t})$,
$\mathbf{t}\in \mathbf{L_{p}}\diagdown \mathbf{L_{p'}}$,
is less than $4-\varepsilon$ for sufficiently small $\varepsilon>0$
and that
$$E[(X(\mathbf{0})+X(\mathbf{t}))-(X(\mathbf{0})+X(\mathbf{s}))]^{2}=2(|t_{1}-s_{1}|^{\alpha_{1}}+|t_{2}-s_{2}|^{\alpha_{1}})(1+o(1))$$
as $\mathbf{t}-\mathbf{s}\rightarrow\mathbf{0}$, so by Theorem 8.1
of Piterbarg (1996), for all sufficiently large $u$ and some
positive $\varepsilon'<\varepsilon$,
\begin{eqnarray*}
P\left(\max_{\mathbf{t}\in \mathbf{L_{p}}\diagdown \mathbf{L_{p'}}}(X(\mathbf{0})+X(\mathbf{t}))>2u+w\right)
&\leq& Cp_{1}p_{2}(2u+w)^{2/\alpha_{1}+2/\alpha_{2}}\Psi\left(\frac{2u+w}{\sqrt{4-\varepsilon}}\right)\\
&\leq& C(u)^{2/\alpha_{1}+2/\alpha_{2}-1}\exp\left(-\frac{u^{2}}{2-\varepsilon'/2}\right)\\
&=&o(\Psi(u)),
\end{eqnarray*}
as $u\rightarrow\infty$.\\
We will apply Theorem 8.2 of Piterbarg (1996), for the first
probability in the right-hand part of (\ref{eq312}). To this end, by
some simple calculations, we get for the correlation function of the
field $X(\mathbf{0})+X(\mathbf{t})$, $\mathbf{t}\in
\mathbf{L_{p'}}$
\begin{eqnarray*}
1-\frac{E(X(\mathbf{0})+X(\mathbf{t}))(X(\mathbf{0})+X(\mathbf{s}))}{\sqrt{E(X(\mathbf{0})+X(\mathbf{t}))^{2}E(X(\mathbf{0})+X(\mathbf{s}))^{2}}}
&\leq& \frac{1-r(\mathbf{t}-\mathbf{s})}{2\sqrt{1+r(\mathbf{t})}\sqrt{1+r(\mathbf{s})}}\\
&\leq& \frac{|t_{1}-s_{1}|^{\alpha_{1}}+|t_{2}-s_{2}|^{\alpha_{2}}}{2(2-\delta_{1}'^{\alpha_{1}}-\delta_{2}'^{\alpha_{2}})}\\
&\leq& 1-\exp(-|t_{1}-s_{1}|^{\alpha_{1}}-|t_{2}-s_{2}|^{\alpha_{2}}),
\end{eqnarray*}
where we assume an additionally that
$p_{1}'^{\alpha_{1}}+p_{2}'^{\alpha_{2}}\leq3/2$. For the
variance of the field $X(\mathbf{0})+X(\mathbf{t})$, $\mathbf{t}\in
\mathbf{L_{p'}}$ we have
$$Var(X(\mathbf{0})+X(\mathbf{t}))=2+2r(\mathbf{t})=4-2(|t_{1}|^{\alpha_{1}}+|t_{2}|^{\alpha_{2}})(1+o(1))$$
as $\mathbf{t}\rightarrow0$, and the point $\mathbf{t}=\mathbf{0}$
is the unique point of maximum of variance of the field
$X(\mathbf{0})+X(\mathbf{t})$. By Slepian's inequality
\begin{eqnarray*}
&&P\left(\max_{\mathbf{t}\in\mathbf{L_{p'}}}(X(\mathbf{0})+X(\mathbf{t}))>2u+w\right)\\
&&=P\left(\max_{\mathbf{t}\in\mathbf{L_{p'}}}\frac{X(\mathbf{0})+X(\mathbf{t})}{\sqrt{E(X(\mathbf{0})+X(\mathbf{t}))^{2}}}
    \sqrt{E(X(\mathbf{0})+X(\mathbf{t}))^{2}}>2u+w\right)\\
&&\leq P\left(\max_{\mathbf{t}\in\mathbf{L_{p'}}}Y(\mathbf{t})
    \sqrt{E(X(\mathbf{0})+X(\mathbf{t}))^{2}}>2u+w\right),
\end{eqnarray*}
where $Y(\mathbf{t})$ is a Gaussian zero mean homogeneous field with
covariance function
$\exp(-[|t_{1}|^{\alpha_{1}}+|t_{2}|^{\alpha_{2}}])$, and thus the
conditions of Theorem 8.2 of Piterbarg (1996), for the case $(ii)$
holds. By this theorem, for some constants $C,C',$
\begin{eqnarray*}
&&P\left(\max_{\mathbf{t}\in\mathbf{L_{p'}}}Y(\mathbf{t}) \sqrt{E(X(\mathbf{0})+X(\mathbf{t}))^{2}}>2u+w\right)\\
&&= C\Psi(u+w/2)(1+o(1))\\
&&= C'u^{-1} \exp(-u^{2}/2-uw)(1+o(1))\\
&&= C'\Psi(u) \exp\left(-\frac{1}{2}\big((\frac{2}{\alpha_{1}}+\frac{2}{\alpha_{2}})\log
u+\log(p_{1}p_{2})+x\big)\right)(1+o(1))\\
&&= C'\Psi(u)e^{-1/2x}\left(u^{\frac{2}{\alpha_{1}}+\frac{2}{\alpha_{2}}}p_{1}p_{2}\right)^{-1/2}(1+o(1))\\
&&= C'\Psi(u)e^{-1/2x}\left((2\log T_{1}T_{2})^{\frac{1}{\alpha_{1}}+\frac{1}{\alpha_{2}}}p_{1}p_{2}\right)^{-1/2}(1+o(1)).
\end{eqnarray*}
Since $(2\log T_{1}T_{2})^{1/\alpha_{i}}p_{i}\rightarrow\infty$ for sparse grids $\mathfrak{R}(p_{i})$, we
get the assertion of the lemma. \hfill$\Box$

Now we consider the probability
$$P_{\mathbf{S}}(u,x)=P\left(\max_{\mathbf{t}\in \mathbf{I_{S}}\cap \mathfrak{R}(p_{1})\times \mathfrak{R}(p_{2})}X(\mathbf{t})>u,
\max_{\mathbf{t}\in
\mathbf{I_{S}}}X(\mathbf{t})>u+\frac{v^{2}+x}{u}\right)$$ when we
will allow $S_{1}S_{2}$ tends to infinity with $u$ but
not too fast. Define
$$\delta(\varepsilon)=\inf_{\max\{|t_{1}|,|t_{2}|\}\geq \varepsilon}(1-r(\mathbf{t})).$$
Note that $\delta(\varepsilon)$ is positive for all positive
$\varepsilon$.

\textbf{Lemma A2}. {\sl Let $S_{i}=S_{i}(u)\geq 2p_{i}$
for all $u$, $i=1,2$ and
$S_{1}S_{2}u^{2/\alpha_{1}+2/\alpha_{2}}=o(\exp(u^{2}\delta(\varepsilon)/8))$
as $u\rightarrow\infty$. Then there exists an $\varepsilon>0$ such
that
\begin{eqnarray}
\label{eq313}
P\left(\max_{\mathbf{t}\in \mathbf{I_{S}}\cap \mathfrak{R}(p_{1})\times\mathfrak{R}(p_{2})}X(\mathbf{t})>u\right)
\thicksim S_{1}S_{2}p_{1}^{-1}p_{2}^{-1}\Psi(u),
\end{eqnarray}
\begin{eqnarray}
\label{eq314}
P\left(\max_{\mathbf{t}\in \mathbf{I_{S}}}X(\mathbf{t})>u+\frac{v^{2}+x}{u}\right)
\thicksim S_{1}S_{2}p_{1}^{-1}p_{2}^{-1}e^{-x}\mathcal{H}_{\alpha_{1}}\mathcal{H}_{\alpha_{2}}\Psi(u),
\end{eqnarray}
as $u\rightarrow\infty$ and
\begin{eqnarray*}
P_{\mathbf{S}}(u,x)=o\left(P\left(\max_{\mathbf{t}\in \mathbf{I_{S}}\cap \mathfrak{R}(p_{1})\times\mathfrak{R}(p_{2})}X(\mathbf{t})>u\right)
+P\left(\max_{\mathbf{t}\in \mathbf{I_{S}}}X(\mathbf{t})>u+\frac{v^{2}+x}{u}\right)\right)
\end{eqnarray*}
as $u\rightarrow\infty$ so that
\begin{eqnarray}
\label{eq315}
&&1-P\left(\max_{\mathbf{t}\in \mathbf{I_{S}}\cap \mathfrak{R}(p_{1})\times\mathfrak{R}(p_{2})}X(\mathbf{t})\leq u,
\max_{\mathbf{t}\in \mathbf{I_{S}}}X(\mathbf{t})\leq u+\frac{v^{2}+x}{u}\right)\nonumber\\
&& \thicksim P\left(\max_{\mathbf{t}\in \mathbf{I_{S}}\cap \mathfrak{R}(p_{1})\times\mathfrak{R}(p_{2})}X(\mathbf{t})>u\right)
    +P\left(\max_{\mathbf{t}\in \mathbf{I_{S}}}X(\mathbf{t})>u+\frac{v^{2}+x}{u}\right)\nonumber\\
&& \thicksim S_{1}S_{2}p_{1}^{-1}p_{2}^{-1}\Psi(u)(1+e^{-x}\mathcal{H}_{\alpha_{1}}\mathcal{H}_{\alpha_{2}})
\end{eqnarray}
as $u\rightarrow\infty$.}

\textbf{Proof:} The relation (\ref{eq314}) is in fact a special case
of Theorem 7.2 of Piterbarg (1996) and relation (\ref{eq313}) can be proved by the same way. Now we prove that for a sparse
grid the double probability $P_{\mathbf{S}}(u,x)$ tends to zero
faster than right-hand part of (\ref{eq314}). Let
$\mathbf{J_{l}}=[(l_{1}-1)p_{1},
(l_{1}+1)p_{1}]\times[(l_{2}-1)p_{2},
(l_{2}+1)p_{2}]$, where $l_{i}=0,1,2,\cdots,
[S_{i}/p_{i}]$, $i=1,2$. We have
\begin{eqnarray}
\label{eq316}
P_{\mathbf{S}}(u,x)&\leq& \sum_{k_{1},l_{1}=0}^{[S_{1}/p_{1}]}\sum_{k_{2},l_{2}=0}^{[S_{2}/p_{2}]}
P\left(X(k_{1}p_{1},k_{2}p_{2})>u, \max_{\mathbf{t}\in \mathbf{J_{l}}}X(\mathbf{t})>u+\frac{v^{2}+x}{u}\right)
=:\sum_{k_{1},l_{1}=0}^{[S_{1}/p_{1}]}\sum_{k_{2},l_{2}=0}^{[S_{2}/p_{2}]}P_{\mathbf{k},\mathbf{l}}\nonumber\\
&=&\sum_{k_{1},l_{1}=0, |k_{1}-l_{1}|\leq 1}^{[S_{1}/p_{1}]}\sum_{k_{2},l_{2}=0, |k_{2}-l_{2}|\leq 1}^{[S_{2}/p_{2}]}P_{\mathbf{k},\mathbf{l}}
+\sum_{k_{1},l_{1}=0, |k_{1}-l_{1}|\leq 1}^{[S_{1}/p_{1}]}\sum_{k_{2},l_{2}=0, |k_{2}-l_{2}|> 1}^{[S_{2}/p_{2}]}P_{\mathbf{k},\mathbf{l}}\nonumber\\
&&+\sum_{k_{1},l_{1}=0, |k_{1}-l_{1}|> 1}^{[S_{1}/p_{1}]}\sum_{k_{2},l_{2}=0, |k_{2}-l_{2}|\leq 1}^{[S_{2}/p_{2}]}P_{\mathbf{k},\mathbf{l}}
+\sum_{k_{1},l_{1}=0, |k_{1}-l_{1}|> 1}^{[S_{1}/p_{1}]}\sum_{k_{2},l_{2}=0, |k_{2}-l_{2}|>1}^{[S_{2}/p_{2}]}P_{\mathbf{k},\mathbf{l}}.
\end{eqnarray}
The members of the first term on the right-hand side of
(\ref{eq316}) can be estimated by Lemma A1, so that
\begin{eqnarray}
\label{eq317}
\sum_{k_{1},l_{1}=0, |k_{1}-l_{1}|\leq 1}^{[S_{1}/p_{1}]}\sum_{k_{2},l_{2}=0, |k_{2}-l_{2}|\leq 1}^{[S_{2}/p_{2}]}P_{\mathbf{k},\mathbf{l}}
=S_{1}S_{2}p_{1}^{-1}p_{2}^{-1}o(\Psi(u))
\end{eqnarray}
as $u\rightarrow\infty$. Let $\mathbf{m}=(m_{1},m_{2})$ with
$m_{i}=0,1,2,\cdots, [S_{i}/p_{i}]$, $i=1,2$.
 Now consider the probability
$P_{\mathbf{k},\mathbf{k+m}}=P_{\mathbf{0},\mathbf{m}}$ for
$\max\{m_{1},m_{2}\}>1$. We have, using Theorem 8.1 of Piterbarg
(1996),
\begin{eqnarray*}
P_{\mathbf{0},\mathbf{m}}&\leq& P\left(\max_{\mathbf{t}\in \mathbf{J_{m}}}(X(\mathbf{0})+X(\mathbf{t}))>2u+\frac{v^{2}+x}{u}\right)\\
&\leq& Cp_{1}p_{2}u^{\frac{2}{\alpha_{1}}+\frac{2}{\alpha_{2}}-1}
\exp\left(-\frac{(2u+(v^{2}+x)/u)^{2}}{2\max_{\mathbf{t}\in \mathbf{J_{m}}}(2+2r(\mathbf{t}))}\right)\\
&\leq& Cp_{1}p_{2}u^{\frac{2}{\alpha_{1}}+\frac{2}{\alpha_{2}}-1}
\exp\left(-\frac{u^{2}+v^{2}}{2(1-\frac{1}{2}\min_{\mathbf{t}\in \mathbf{J_{m}}}(1-r(\mathbf{t}))}\right)\\
&\leq& Cp_{1}p_{2}u^{\frac{2}{\alpha_{1}}+\frac{2}{\alpha_{2}}-1}
\exp\left(-\frac{1}{2}(u^{2}+v^{2})(1+\frac{1}{2}\min_{\mathbf{t}\in \mathbf{J_{m}}}(1-r(\mathbf{t}))\right)\\
&\leq& Cp_{1}^{1/2}p_{2}^{1/2}u^{\frac{1}{\alpha_{1}}+\frac{1}{\alpha_{2}}-1}
\exp\left(-\frac{1}{2}u^{2}\right)\exp\left(-\frac{1}{4}u^{2}\min_{\mathbf{t}\in \mathbf{J_{m}}}(1-r(\mathbf{t}))\right)\\
&\leq& Cp_{1}^{1/2}p_{2}^{1/2}u^{\frac{1}{\alpha_{1}}+\frac{1}{\alpha_{2}}}
\Psi(u)\exp\left(-\frac{1}{4}u^{2}\min_{\mathbf{t}\in \mathbf{J_{m}}}(1-r(\mathbf{t}))\right).
\end{eqnarray*}
Let $\varepsilon$ be such that $1-r(\mathbf{t})\geq
\frac{1}{2}(|t_{1}|^{\alpha_{1}}+|t_{2}|^{\alpha_{2}})$ for all
$\mathbf{t}\in
(-\varepsilon,\varepsilon)\times(-\varepsilon,\varepsilon)$. Then
$$P_{\mathbf{0},\mathbf{m}}\leq Cp_{1}^{1/2}p_{2}^{1/2}u^{\frac{1}{\alpha_{1}}+\frac{1}{\alpha_{2}}}
\Psi(u)\exp\left(-\frac{1}{8}u^{2}\delta(\varepsilon)\right)$$
for $\max\{|(m_{1}-1)p_{1}|, |(m_{2}-1)p_{2}|\}>\varepsilon$ and
$$P_{\mathbf{0},\mathbf{m}}\leq Cp_{1}^{1/2}p_{2}^{1/2}u^{\frac{1}{\alpha_{1}}+\frac{1}{\alpha_{2}}}
\Psi(u)\exp\left(-\frac{1}{8}u^{2}[|(m_{1}-1)p_{1}|^{\alpha_{1}}+|(m_{2}-1)p_{2}|^{\alpha_{2}}]\right)$$
for $\max\{|(m_{1}-1)p_{1}|, |(m_{2}-1)p_{2}|\}\leq\varepsilon$.
Thus, letting $\mathbf{i}=\mathbf{l}-\mathbf{k}$, for the second sum we
have
\begin{eqnarray*}
\sum_{k_{1},l_{1}=0, |k_{1}-l_{1}|\leq 1}^{[S_{1}/p_{1}]}\sum_{k_{2},l_{2}=0, |k_{2}-l_{2}|> 1}^{[S_{2}/p_{2}]}P_{\mathbf{k},\mathbf{l}}
&\leq&4S_{1}S_{2}p_{1}^{-1}p_{2}^{-1}\sum_{i_{1}=0}^{1}\sum_{i_{2}=2}^{[S_{2}/p_{2}]}p_{\mathbf{0},\mathbf{i}}\\
&\leq& CS_{1}S_{2}p_{1}^{-1/2}p_{2}^{-1/2}u^{\frac{1}{\alpha_{1}}+\frac{1}{\alpha_{2}}}\Psi(u)
\bigg\{S_{2}p_{2}^{-1}\exp\left(-\frac{1}{8}u^{2}\delta(\varepsilon)\right)\\
&&+\sum_{i_{1}-1=-1}^{0}\sum_{i_{2}-1=1}^{\varepsilon/p_{2}}
\exp\left(-\frac{1}{8}u^{2}[|(i_{1}-1)p_{1}|^{\alpha_{1}}+|(i_{2}-1)p_{2}|^{\alpha_{2}}]\right)\bigg\}\\
&\leq&CS_{1}S_{2}p_{1}^{-1}p_{2}^{-1}\Psi(u)o(1)
\end{eqnarray*}
as $u\rightarrow\infty$. Similarly,
$$\sum_{k_{1},l_{1}=0, |k_{1}-l_{1}|> 1}^{[S_{1}/p_{1}]}\sum_{k_{2},l_{2}=0, |k_{2}-l_{2}|\leq 1}^{[S_{2}/p_{2}]}P_{\mathbf{k},\mathbf{l}}
\leq CS_{1}S_{2}p_{1}^{-1}p_{2}^{-1}\Psi(u)o(1)$$ as
$u\rightarrow\infty$. For the fourth sum, we have
\begin{eqnarray*}
\sum_{k_{1},l_{1}=0, |k_{1}-l_{1}|> 1}^{[S_{1}/p_{1}]}\sum_{k_{2},l_{2}=0, |k_{2}-l_{2}|>1}^{[S_{2}/p_{2}]}P_{\mathbf{k},\mathbf{l}}
&\leq&4S_{1}S_{2}p_{1}^{-1}p_{2}^{-1}\sum_{i_{1}=2}^{[S_{1}/\delta_{1}]}\sum_{i_{2}=2}^{[S_{2}/p_{2}]}P_{\mathbf{0},\mathbf{i}}\\
&\leq& CS_{1}S_{2}p_{1}^{-1/2}p_{2}^{-1/2}u^{\frac{1}{\alpha_{1}}+\frac{1}{\alpha_{2}}}\Psi(u)
\bigg\{S_{1}p_{1}^{-1}S_{2}p_{2}^{-1}\exp\left(-\frac{1}{8}u^{2}\delta(\varepsilon)\right)\\
&&+\sum_{i_{1}-1=1}^{\varepsilon/p_{1}}\sum_{i_{2}-1=1}^{\varepsilon/p_{2}}
\exp\left(-\frac{1}{8}u^{2}[|(i_{1}-1)p_{1}|^{\alpha_{1}}+|(i_{2}-1)p_{2}|^{\alpha_{2}}]\right)\bigg\}\\
&\leq&CS_{1}S_{2}p_{1}^{-1}p_{2}^{-1}\Psi(u)o(1),
\end{eqnarray*}
as $u\rightarrow\infty$. Now we can easily prove the relation
(\ref{eq314}). We have for all $\mathbf{k}$ and $\mathbf{l}$
$$P\left(X(k_{1}p_{1},k_{2}p_{2})>u,X(l_{1}p_{1},l_{2}p_{2})>u\right)\leq P_{\mathbf{k},\mathbf{l}},$$
hence
$$\sum_{k_{1},l_{1}=0,k_{1}\neq l_{1}}^{[S_{1}/p_{1}]}\sum_{k_{2},l_{2}=0,k_{2}\neq l_{2}}^{[S_{2}/p_{2}]}
P\left(X(k_{1}p_{1},k_{2}p_{2})>u,X(l_{1}p_{1},l_{2}p_{2})>u\right)
\leq\sum_{k_{1},l_{1}=0,k_{1}\neq
l_{1}}^{[S_{1}/p_{1}]}\sum_{k_{2},l_{2}=0,k_{2}\neq
l_{2}}^{[S_{2}/p_{2}]}P_{\mathbf{k},\mathbf{l}},$$ from which
it follows that double sum in the above left-hand side tends to zero
faster than $S_{1}S_{2}p_{1}^{-1}p_{2}^{-1}\Psi(u)$ as
$u\rightarrow\infty$. Thus, both the assertions of Lemma A2 are
proved.\hfill$\Box$

\subsection{Pickands grid}

Lat $\mathbf{a}=(a_{1},a_{2})>(0,0)$. In this subsection suppose that
$\mathfrak{R}(p_{i})$, $i=1,2$ are Pickands grids, ie.,
$\mathfrak{R}(p_{i})=\{a_{i}ku^{-2/\alpha_{i}}, k\in
\mathbb{N}\}$. We will evaluate the asymptotic behavior of the
probability
$$P_{\mathbf{S}}'(u,x)=P\left(\max_{\mathbf{t}\in \mathbf{I_{S}}\cap \mathfrak{R}(p_{1})\times \mathfrak{R}(p_{2})}X(\mathbf{t})>u,
\max_{\mathbf{t}\in
\mathbf{I_{S}}}X(\mathbf{t})>u+\frac{x}{u}\right).$$ As in the
previous subsection, we begin with a short interval. Let
$\lambda_{i}>a_{i}$. Then it can be proved quite similar to the
proof of Lemma 6.1 of Piterbarg (1996), that
$$P_{(\lambda_{1}u^{-2/\alpha_{1}}, \lambda_{2}u^{-2/\alpha_{2}})}'(u,x)\thicksim \mathcal{H}_{\mathbf{a},\alpha_{1},\alpha_{2}}^{0,x}\Psi(u)$$
as $u\rightarrow\infty$, where
$$H_{\mathbf{d},\alpha_{1},\alpha_{2}}^{0,x}(\lambda_{1},\lambda_{2})=\int_{-\infty}^{+\infty}e^{s}
P\left(\max_{(k_{1}d_{1},k_{2}d_{2})\in[0,\lambda_{1}]\times[0,\lambda_{2}]}\sqrt{2}\chi(k_{1}d_{1},k_{2}d_{2})>s,
\max_{(t_{1},t_{2})\in[0,\lambda_{1}]\times[0,\lambda_{2}]}\sqrt{2}\chi(t_{1},t_{2})>s+x\right)ds.$$
It also can be proved in a similar way as for Lemma 6.1 and Theorem 7.2 of Piterbarg
(1996) that
$$\mathcal{H}_{\mathbf{a},\alpha_{1},\alpha_{2}}^{0,x}:=\lim_{\lambda_{1}\rightarrow\infty\atop\lambda_{2}\rightarrow\infty}
\mathcal{H}_{\mathbf{a},\alpha_{1},\alpha_{2}}^{0,x}(\lambda_{1},\lambda_{2})/(\lambda_{1}\lambda_{2})\in
(0,\infty)$$ and that there exsits $\kappa\in (0,1/2)$ such that
for any $S_{i}=S_{i}(u)$ with
$S_{1}S_{2}u^{-(2/\alpha_{1}+2/\alpha_{2})}\rightarrow\infty$ and
$S_{1}S_{2}=O(\exp(\kappa u^{2}))$ as
$u\rightarrow\infty$ with
\begin{eqnarray}
\label{eq318}
P_{\mathbf{S}}'(u,x)\thicksim S_{1}S_{2}\mathcal{H}_{\mathbf{a},\alpha_{1},\alpha_{2}}^{0,x}u^{2/\alpha_{1}+2/\alpha_{2}}\Psi(u)
\end{eqnarray}
as $u\rightarrow\infty$, respectively.  From here we have for Pickands grids,

\textbf{Lemma A3}. {\sl For any $\mathbf{a}=(a_{1},a_{2})$ and
$\mathfrak{R}(p_{i})=\{a_{i}ku^{-2/\alpha_{i}}, k\in
\mathbb{N}\}$,
\begin{eqnarray}
\label{eq319}
&&1-P\left(\max_{\mathbf{t}\in \mathbf{I_{S}}\cap \mathfrak{R}(p_{1})\times \mathfrak{R}(p_{2})}X(\mathbf{t})\leq u,
\max_{\mathbf{t}\in \mathbf{I_{S}}}X(\mathbf{t})\leq u+\frac{v^{2}+x}{u}\right)\nonumber\\
&&=P\left(\max_{\mathbf{t}\in \mathbf{I_{S}}\cap \mathfrak{R}(p_{1})\times \mathfrak{R}(p_{2})}X(\mathbf{t})> u\right)
+P\left(\max_{\mathbf{t}\in \mathbf{I_{S}}}X(\mathbf{t})> u+\frac{v^{2}+x}{u}\right)\nonumber\\
&&\ \ \ \ -P\left(\max_{\mathbf{t}\in \mathbf{I_{S}}\cap \mathfrak{R}(p_{1})\times \mathfrak{R}(p_{2})}X(\mathbf{t})>u,
\max_{\mathbf{t}\in\mathbf{I_{S}}}X(\mathbf{t})>u+\frac{x}{u}\right)\nonumber\\
&&\thicksim S_{1}S_{2}\mathcal{H}_{\alpha_{1}}\mathcal{H}_{\alpha_{2}}\left(u+\frac{v^{2}+x}{u}\right)^{2/\alpha_{1}+2/\alpha_{2}}\Psi\left(u+\frac{v^{2}+x}{u}\right)
+S_{1}S_{2}\mathcal{H}_{a_{1},\alpha_{1}}\mathcal{H}_{a_{2},\alpha_{2}}u^{2/\alpha_{1}+2/\alpha_{2}}\Psi(u)\nonumber\\
&&\ \ \ \ -S_{1}S_{2}\mathcal{H}_{\mathbf{a},\alpha_{1},\alpha_{2}}^{0,x}u^{2/\alpha_{1}+2/\alpha_{2}}\Psi(u)
\end{eqnarray}
as $u\rightarrow\infty$. }

\subsection{Dense grid}

In this subsection, we state a lemma for the dense grid case which is important for our proofs.

\textbf{Lemma A4}. {\sl Let $S_{i}=S_{i}(u)$ with
$S_{1}S_{2}u^{-(2/\alpha_{1}+2/\alpha_{2})}\rightarrow\infty$ and
$S_{1}S_{2}=O(\exp(\kappa u^{2}))$ with $\kappa\in (0,1/2]$ as
$u\rightarrow\infty$.  For any $\mathbf{a}=(a_{1},a_{2})$ and
$\mathfrak{R}(p_{i})=\{a_{i}ku^{-2/\alpha_{i}}, k\in
\mathbb{N}\}$, we have
\begin{eqnarray}
\label{eq320}
P\left(\max_{\mathbf{t}\in \mathbf{I_{S}}\cap \mathfrak{R}(p_{1})\times \mathfrak{R}(p_{2})}X(\mathbf{t})\leq u\right)
-P\left(\max_{\mathbf{t}\in \mathbf{I_{S}}}X(\mathbf{t})\leq u\right)
= g(a_{1},a_{2})\mathcal{H}_{\alpha_{1}}\mathcal{H}_{\alpha_{2}}S_{1}S_{2}u^{2/\alpha_{1}+2/\alpha_{2}}\Psi(u),
\end{eqnarray}
where $g(a_{1},a_{2})\rightarrow 0$ as $\mathbf{a}\rightarrow \mathbf{0}$.}

\textbf{Proof:} Lemma 1 of D\c{e}bicki et
al. (2014) shows that (\ref{eq320}) holds for some fixed $S_{i}>0$. By the homogeneity of $X(\mathbf{t})$, it is easy to
extend  (\ref{eq320}) to the case $S_{1}S_{2}=O(\exp(\kappa u^{2}))$ with $\kappa\in (0,1/2]$, see eg. the proof of Lemma 12.3.2
of Leadbetter (1983) for more details.   \hfill$\Box$

\section{Appendix B}

In this section, we give three technical lemmas which are used for the proof of Lemma 3.1. Recall that $u_{\mathbf{T}}=
b_{\mathbf{T}}+x/a_{\mathbf{T}}$, $u_{\mathbf{T}}'=
b_{\mathbf{T}}'+y/a_{\mathbf{T}}$, where
$b_{\mathbf{T}}^{'}=b_{\mathbf{T}}^{\mathbf{p}}$ for  sparse grids and
$b_{\mathbf{T}}^{'}=b_{\mathbf{a},\mathbf{T}}$ for  Pickands grids, and $r^{(h)}(\mathbf{kq},\mathbf{lq})=hr(\mathbf{kq},\mathbf{lq})+(1-h)\varrho(\mathbf{kq},\mathbf{lq})$ with $h\in[0,1]$.
Let
$$\varpi(\mathbf{t},\mathbf{s})=\max\{|r(\mathbf{t},\mathbf{s})|, |\varrho(\mathbf{t},\mathbf{s})|\}$$
and
 $$\vartheta(\mathbf{\mathbf{z}})=\sup_{\mathbf{0}\leq\mathbf{s,t}\leq \mathbf{T},\atop \{|s_{1}-t_{1}|>z_{1}\}\cup\{|s_{2}-t_{2}|>z_{2}\}}\{\varpi(\mathbf{t},\mathbf{s})\}.
$$
It is easy to see from Assumptions {\bf A1} and {\bf A2} that for any $\varepsilon_{1}>0$ and $\varepsilon_{2}>0$
$$\vartheta(\varepsilon_{1},\varepsilon_{2})<1$$
for all sufficiently large $\mathbf{T}$. Further, let $a,b$ be such that
$$0<b<a<\big(1-\vartheta(\varepsilon,\varepsilon)\big)/\big(1+\vartheta(\varepsilon,\varepsilon)\big)<1$$
 for all sufficiently large $\mathbf{T}$ and for some $\varepsilon>0$ which will be chosen in the blow.\\

\textbf{Lemma B1}. {\sl Under the conditions of Lemma 3.3, we have
\begin{eqnarray}
\label{eqB11}
\sum_{\mathbf{kq}\in \mathbf{O_{i}},\mathbf{lq}\in \mathbf{O_{j}}\atop \mathbf{kq}\neq \mathbf{lq}, \mathbf{1}\leq \mathbf{i,j}\leq \mathbf{n} }|r(\mathbf{kq},\mathbf{lq})-\varrho(\mathbf{kq},\mathbf{lq})|
\int_{0}^{1}\frac{1}{\sqrt{1-r^{(h)}(\mathbf{kq},\mathbf{lq})}}\exp\left(-\frac{u_{\mathbf{T}}^{2}}{1+r^{(h)}(\mathbf{kq},\mathbf{lq})}\right)dh
\rightarrow 0
\end{eqnarray}
as $\mathbf{T}\rightarrow\infty$.
}

\textbf{Proof:} Recall that  $\mathfrak{R}(q_{i})$, $i=1,2$ are Pickands grids.
First, we consider the case that $\mathbf{kq},\mathbf{lq}$ in the
same interval $\mathbf{O_{i}}$.
Split the sum (\ref{eqB11}) into two parts as
\begin{eqnarray}
\label{eq402}
\sum_{\mathbf{kq},\mathbf{lq}\in \mathbf{O_{i}},\mathbf{kq}\neq\mathbf{lq},\mathbf{i}=\mathbf{1},\cdots,\mathbf{n},\atop \max\{|l_{1}q_{1}-k_{1}q_{1}|,|l_{2}q_{2}-k_{2}q_{2}|\}\leq\varepsilon}
+\sum_{\mathbf{kq},\mathbf{lq}\in \mathbf{O_{i}},\mathbf{kq}\neq\mathbf{lq}, \mathbf{i}=\mathbf{1},\cdots,\mathbf{n},\atop \max\{|l_{1}q_{1}-k_{1}q_{1}|,|l_{2}q_{2}-k_{2}q_{2}|\}>\varepsilon}=:J_{\mathbf{T},1}+J_{\mathbf{T},2}.
\end{eqnarray}
We deal with $J_{\mathbf{T},1}$ and note that in this case, by the definition of the field $\xi_{\mathbf{T}}(\mathbf{t})$, we have
$\varrho(\mathbf{kq},\mathbf{lq})-r(\mathbf{kq},\mathbf{lq})=\rho(\mathbf{T})(1-r(\mathbf{kq},\mathbf{lq}))$.
By  Assumption {\bf A1}  we can choose small enough $\varepsilon>0$ such that
$\varrho(\mathbf{kq},\mathbf{lq})=r(\mathbf{kq},\mathbf{lq})+(1-r(\mathbf{kq},\mathbf{lq}))\rho(\mathbf{T})\sim
r(\mathbf{kq},\mathbf{lq})$ for sufficiently large $\mathbf{T}$ and $\max\{|l_{1}q_{1}-k_{1}q_{1}|,|l_{2}q_{2}-k_{2}q_{2}|\}\leq\varepsilon$.
It follows from Assumption \textbf{A1} again that for all
$|t_{i}|\leq\varepsilon<2^{-1/\alpha_{i}}$,
\begin{eqnarray}
\label{eqT20}
\frac{1}{2}(|t_{1}|^{\alpha_{1}}+|t_{2}|^{\alpha_{2}})\leq 1-r(\mathbf{t})\leq 2(|t_{1}|^{\alpha_{1}}+|t_{2}|^{\alpha_{2}})
\end{eqnarray}
and the definition of $u_{\mathbf{T}}$ implies
\begin{eqnarray}
\label{eq403}
u^{2}_{\mathbf{T}}=2\log T_{1}T_{2}-\log\log T_{1}T_{2}+(\frac{2}{\alpha_{1}}+\frac{2}{\alpha_{2}})\log\log T_{1}T_{2}+O(1).
\end{eqnarray}
Consequently, since further $q_{i}=\gamma_{i}(\log
T_{1}T_{2})^{-1/\alpha_{i}}$ we obtain
%(without loss of generality we
%assume that $\mathbf{k}\geq \mathbf{0}$ and $\mathbf{k}\neq\mathbf{0}$)
\begin{eqnarray}
\label{eq404}
J_{\mathbf{T},1}&\leq&C\sum_{\mathbf{kq},\mathbf{lq}\in \mathbf{O_{i}},\mathbf{kq}\neq\mathbf{lq},\mathbf{i}=\mathbf{1},\cdots,\mathbf{n},\atop \max\{|l_{1}q_{1}-k_{1}q_{1}|,|l_{2}q_{2}-k_{2}q_{2}|\}\leq\varepsilon}
|r(\mathbf{kq},\mathbf{lq})-\varrho(\mathbf{kq},\mathbf{lq})|\frac{1}{\sqrt{1-r(\mathbf{kq},\mathbf{lq})}}\exp\left(-\frac{u^{2}_{\mathbf{T}}}{1+r(\mathbf{kq},\mathbf{lq})}\right)\nonumber\\
&\leq & C\frac{T_{1}}{q_{1}}\frac{T_{2}}{q_{2}}\rho(\mathbf{T})\sum_{0<k_{1}q_{1}\leq \varepsilon, 0<k_{2}q_{2}\leq\varepsilon}|1-r(\mathbf{kq})|
\frac{1}{\sqrt{1-r(\mathbf{kq})}}\exp\left(-\frac{u^{2}_{\mathbf{T}}}{2}\right)\exp\left(-\frac{(1-r(\mathbf{kq}))u^{2}_{\mathbf{T}}}{2(1+r(\mathbf{kq}))}\right)\nonumber\\
&\leq & C \frac{T_{1}}{q_{1}}\frac{T_{2}}{q_{2}}\rho(\mathbf{T})T_{1}^{-1}T_{2}^{-1}(\log T_{1}T_{2})^{1/2-1/\alpha_{1}-1/\alpha_{2}}
\sum_{0<k_{1}q_{1}\leq \varepsilon, 0<k_{2}q_{2}\leq\varepsilon}\sqrt{1-r(\mathbf{kq})}\exp\left(-\frac{(1-r(\mathbf{kq}))u^{2}_{\mathbf{T}}}{2(1+r(\mathbf{kq}))}\right)\nonumber\\
&\leq & C (\log T_{1}T_{2})^{-1/2}\sum_{0<k_{1}q_{1}\leq \varepsilon, 0<k_{2}q_{2}\leq\varepsilon}[(kq_{1})^{\alpha_{1}}+(kq_{2})^{\alpha_{2}}]^{1/2}\exp\left(-\frac{1}{4}[(kq_{1})^{\alpha_{1}}+(kq_{2})^{\alpha_{2}}]\log (T_{1}T_{2})\right)\nonumber\\
&\leq & C (\log T_{1}T_{2})^{-1/2}\sum_{0<k_{1}q_{1}\leq \varepsilon, 0<k_{2}q_{2}\leq\varepsilon}\exp\left(-\frac{1}{4}[(kq_{1})^{\alpha_{1}}+(kq_{2})^{\alpha_{2}}]\log (T_{1}T_{2})\right)\nonumber\\
&\leq & C (\log T_{1}T_{2})^{-1/2}\sum_{k_{1}=1}^{\infty}e^{-\frac{1}{4}(k_{1}\gamma_{1})^{\alpha_{1}}}
\sum_{k_{2}=1}^{\infty}e^{-\frac{1}{4}(k_{2}\gamma_{2})^{\alpha_{2}}}\nonumber\\
&\leq & C (\log T_{1}T_{2})^{-1/2},
\end{eqnarray}
which shows $J_{\mathbf{T},1}\rightarrow0$ as $\mathbf{T}\rightarrow\infty$.\\
Using the fact that $u_{\mathbf{T}}\thicksim (2\log
T_{1}T_{2})^{1/2}$, we obtain
\begin{eqnarray}
\label{eq405}
J_{\mathbf{T},2}&\leq&C\sum_{\mathbf{kq},\mathbf{lq}\in \mathbf{O_{i}},\mathbf{kq}\neq\mathbf{lq}, \mathbf{i}=\mathbf{1},\cdots,\mathbf{n},\atop \max\{|l_{1}q_{1}-k_{1}q_{1}|,|l_{2}q_{2}-k_{2}q_{2}|\}>\varepsilon}
      |r(\mathbf{kq},\mathbf{lq})-\varrho(\mathbf{kq},\mathbf{lq})|\exp\left(-\frac{u^{2}_{\mathbf{T}}}{1+\varpi(\mathbf{kq},\mathbf{lq})}\right)\nonumber\\
&\leq & C\frac{T_{1}}{q_{1}}\frac{T_{2}}{q_{2}}\sum_{0\leq k_{1}q_{1}\leq T_{1}^{a}, 0\leq k_{2}q_{2}\leq T_{2}^{a},\mathbf{i}=\mathbf{1},\cdots,\mathbf{n},\atop \max\{k_{1}q_{1},k_{2}q_{2}\}>\varepsilon}
     \exp\left(-\frac{u^{2}_{\mathbf{T}}}{1+\vartheta(\varepsilon,\varepsilon)}\right)\nonumber\\
&\leq & C\frac{T_{1}}{q_{1}}\frac{T_{2}}{q_{2}}\exp\left(-\frac{u^{2}_{\mathbf{T}}}{1+\vartheta(\varepsilon,\varepsilon)}\right)
     \sum_{0\leq k_{1}q_{1}\leq T_{1}^{a}, 0\leq k_{2}q_{2}\leq T_{2}^{a}}1\nonumber\\
&\leq & C \frac{T_{1}}{q_{1}}\frac{T_{2}}{q_{2}}(T_{1}T_{2})^{-\frac{2}{1+\vartheta(\varepsilon,\varepsilon)}}
    \sum_{0\leq k_{1}q_{1}\leq T_{1}^{a}, 0\leq k_{2}q_{2}\leq T_{2}^{a}}1\nonumber\\
&\leq & C (T_{1}T_{2})^{a-\frac{1-\vartheta(\varepsilon,\varepsilon)}{1+\vartheta(\varepsilon,\varepsilon)}}(\log T_{1}T_{2})^{2/\alpha_{1}+2/\alpha_{2}}.
\end{eqnarray}
Thus, $J_{\mathbf{T},2}\rightarrow0$ as $\mathbf{T}\rightarrow\infty$ since $a<\frac{1-\vartheta(\varepsilon,\varepsilon)}{1+\vartheta(\varepsilon,\varepsilon)}$.

Second, we deal with the case that $\mathbf{kq}\in \mathbf{O_{i}}$
and $\mathbf{lq}\in \mathbf{O_{j}}$, $\mathbf{i}\neq \mathbf{j}$.
Note that in this case, the distance between the points in any two
rectangles $\mathbf{O_{i}}$ and $\mathbf{O_{j}}$ is large than
$T_{1}^{b}$ or $T_{2}^{b}$ and $\varrho(\mathbf{kq},\mathbf{lq})=\rho(\mathbf{T})$ for  $\mathbf{kq}\in \mathbf{O_{i}}$
and $\mathbf{lq}\in \mathbf{O_{j}}$, $\mathbf{i}\neq \mathbf{j}$.

Obviously, the sum in (\ref{eqB11}) is smaller than
\begin{eqnarray}
\label{Tan1}
C\sum_{\mathbf{kq}\in \mathbf{O_{i}},\mathbf{lq}\in \mathbf{O_{j}}\atop \mathbf{kq}\neq \mathbf{lq}, \mathbf{1}\leq \mathbf{i}\neq \mathbf{j}\leq \mathbf{n} }
|r(\mathbf{kq},\mathbf{lq})-\rho(\mathbf{T})|\exp\left(-\frac{u^{2}_{\mathbf{T}}}{1+\varpi(\mathbf{kq},\mathbf{lq})}\right).
%&\leq&C\sum_{\mathbf{kq}\in \mathbf{O_{i}},\mathbf{lq}\in \mathbf{O_{j}}\atop \mathbf{kq}\neq \mathbf{lq}, \mathbf{1}\leq \mathbf{i}\neq \mathbf{j}\leq \mathbf{n} }
%|r(\mathbf{kq},\mathbf{lq})-\rho(\mathbf{T})|\exp\left(-\frac{u^{2}_{\mathbf{T}}}{1+\vartheta(T_{1}^{b},T_{2}^{b})}\right)\nonumber\\
%&\leq&  C\frac{T_{1}}{q_{1}}\frac{T_{2}}{q_{2}}\sum_{\mathbf{0}\leq\mathbf{kq}\leq \mathbf{T}, \mathbf{kq}\neq \mathbf{0}}
%|r(\mathbf{kq})-\rho(\mathbf{T})|\exp\left(-\frac{u^{2}_{\mathbf{T}}}{1+\vartheta(T_{1}^{b},T_{2}^{b})}\right)
\end{eqnarray}
Split the sum of (\ref{Tan1}) into three
parts, the first for $|k_{1}q_{1}-l_{1}q_{1}|>0$ and $|k_{2}q_{2}-l_{2}q_{2}|>0$, the second for
$k_{1}q_{1}-l_{1}q_{1}=0$ and $|k_{2}q_{2}-l_{2}q_{2}|>0$, the third for $k_{2}q_{2}-l_{2}q_{2}=0$ and $|k_{1}q_{1}-l_{1}q_{1}|>0$ and
denote them by $S_{\mathbf{T},i}$, $i=1,2,3,$ respectively. Let $\beta$ be such that
$0<b<a<\beta<\frac{1-\vartheta(\varepsilon,\varepsilon)}{1+\vartheta(\varepsilon,\varepsilon)}$
for all sufficiently large $\mathbf{T}$.

We consider the term $S_{\mathbf{T},1}$ and
split it into two parts as
\begin{eqnarray*}
S_{\mathbf{T},1}= C\sum_{\mathbf{kq}\in \mathbf{O_{i}},\mathbf{lq}\in \mathbf{O_{j}}, \mathbf{kq}\neq \mathbf{lq}, \mathbf{1}\leq \mathbf{i}\neq \mathbf{j}\leq \mathbf{n}\atop |k_{1}q_{1}-l_{1}q_{1}||k_{2}q_{2}-l_{2}q_{2}|\leq (T_{1}T_{2})^{\beta}}
+ C\sum_{\mathbf{kq}\in \mathbf{O_{i}},\mathbf{lq}\in \mathbf{O_{j}},\mathbf{kq}\neq \mathbf{lq}, \mathbf{1}\leq \mathbf{i}\neq \mathbf{j}\leq \mathbf{n}\atop |k_{1}q_{1}-l_{1}q_{1}||k_{2}q_{2}-l_{2}q_{2}|> (T_{1}T_{2})^{\beta}}=: S_{\mathbf{T},11}+S_{\mathbf{T},12}.
\end{eqnarray*}
For $S_{\mathbf{T},11}$, with the similar derivation as for
(\ref{eq405}), we have
\begin{eqnarray}
\label{eq407}
S_{\mathbf{T},11}
&\leq & C\frac{T_{1}}{q_{1}}\frac{T_{2}}{q_{2}}\sum_{0\leq k_{1}q_{1}\leq T_{1}, 0\leq k_{2}q_{2}\leq T_{2},\atop k_{1}q_{1}k_{2}q_{2}\leq (T_{1}T_{2})^{\beta}}
     \exp\left(-\frac{u^{2}_{\mathbf{T}}}{1+\vartheta(\varepsilon,\varepsilon)}\right)\nonumber\\
&\leq & C (T_{1}T_{2})^{\beta-\frac{1-\vartheta(\varepsilon,\varepsilon)}{1+\vartheta(\varepsilon,\varepsilon)}}(\log T_{1}T_{2})^{2/\alpha_{1}+2/\alpha_{2}}.
\end{eqnarray}
Consequently, since $\beta<\frac{1-\vartheta(\varepsilon,\varepsilon)}{1+\vartheta(\varepsilon,\varepsilon)}$, we have $S_{\mathbf{T},11}\rightarrow0$ as $\mathbf{T}\rightarrow\infty$.\\
For $S_{\mathbf{T},12}$, we need more precise estimation. Let's define
$$\omega_{1}(\mathbf{t})=\max\{|r(\mathbf{t})|, |\rho(\mathbf{T})|\}$$
and
$$
\theta_{1}(\mathbf{\mathbf{z}})=\sup_{\mathbf{0}\leq\mathbf{t}\leq \mathbf{T},\atop |t_{1}t_{2}|>z_{1}z_{2}}\{\omega_{1}(\mathbf{t})\}.$$
By the Assumption \textbf{A3}, there exist constants $C>0$ and $K>0$
such that
$$\theta_{1}(\mathbf{t})\log\left(t_{1}t_{2}\right)\leqslant K$$
for all $\mathbf{T}$ sufficiently large and $\mathbf{t}$ satisfying
$t_{1}t_{2}\geq C$. Thus for
all $\mathbf{T}$ large enough and for $(k_{1}q_{1},k_{2}q_{2})$ such
that $k_{1}q_{1}k_{2}q_{2}\geq
(T_{1}T_{2})^\beta$,
$\theta_{1}(\mathbf{kq})\leq K/\log(T_{1}T_{2})^{\beta}$.
Now making use of (\ref{eq403}), we
obtain
\begin{eqnarray}
\label{Tan2}
&&\frac{(T_{1}T_{2})^{2}}{q_{1}^{2}q_{2}^{2}\log (T_{1}T_{2})}\exp\left(-\frac{u_{\mathbf{T}}^{2}}{1+\theta_{1}(T_{1}^{\beta},T_{2}^{\beta})}\right)\nonumber\\
&&\leq \frac{(T_{1}T_{2})^{2}}{q_{1}^{2}q_{2}^{2}\log (T_{1}T_{2})}\exp\left(-\frac{u_{\mathbf{T}}^{2}}{1+K/\log(T_{1}T_{2})^{\beta}}\right)\nonumber\\
&&\thicksim \frac{(T_{1}T_{2})^{2}}{q_{1}^{2}q_{2}^{2}\log (T_{1}T_{2})}\left((T_{1}T_{2})^{-2}(\log T_{1}T_{2}) (\log T_{1}T_{2})^{-(2/\alpha_{1}+2/\alpha_{2})}\right)^{\frac{1}{1+K/\log(T_{1}T_{2})^{\beta}}}\nonumber\\
&&\leq O(1)(T_{1}T_{2})^{(2K/\log(T_{1}T_{2})^{\beta})/(1+K/\log(T_{1}T_{2})^{\beta})}(\log T_{1}T_{2})^{((2/\alpha_{1}+2/\alpha_{2}-1)K/\log (T_{1}T_{2})^{\beta})/(1+K/\log (T_{1}T_{2})^{\beta})}\nonumber\\
&&=O(1).
\end{eqnarray}
Therefore, by a similar argument as for the proof of Lemma 6.4.1 of
Leadbetter et al. (1983) we obtain
\begin{eqnarray}
\label{eq409}
S_{\mathbf{T},12}&\leq&C\frac{T_{1}}{q_{1}}\frac{T_{2}}{q_{2}}
      \sum_{\mathbf{0}\leq\mathbf{kq}\leq \mathbf{T}, \mathbf{kq}\neq \mathbf{0}\atop k_{1}q_{1}k_{2}q_{2}> (T_{1}T_{2})^{\beta}}
      |r(\mathbf{kq})-\rho(\mathbf{T})|\exp\left(-\frac{u^{2}_{\mathbf{T}}}{1+\theta_{1}(T_{1}^{\beta},T_{2}^{\beta})}\right)\nonumber\\
&\leq&C\frac{T_{1}}{q_{1}}\frac{T_{2}}{q_{2}}\exp\left(-\frac{u^{2}_{\mathbf{T}}}{1+\theta_{1}(T_{1}^{\beta},T_{2}^{\beta})}\right)
      \sum_{\mathbf{0}\leq\mathbf{kq}\leq \mathbf{T}, \mathbf{kq}\neq \mathbf{0}\atop k_{1}q_{1}k_{2}q_{2}> (T_{1}T_{2})^{\beta}}|r(\mathbf{kq})-\rho(\mathbf{T})|\nonumber\\
&=& C\frac{(T_{1}T_{2})^{2}}{q_{1}^{2}q_{2}^{2}\log (T_{1}T_{2})}\exp\left(-\frac{u_{\mathbf{T}}^{2}}{1+\theta_{1}(T_{1}^{\beta},T_{2}^{\beta})}\right)\cdot \frac{q_{1}q_{2}\log (T_{1}T_{2})}{T_{1}T_{2}}
      \sum_{\mathbf{0}\leq\mathbf{kq}\leq \mathbf{T}, \mathbf{kq}\neq \mathbf{0}\atop k_{1}q_{1}k_{2}q_{2}> (T_{1}T_{2})^{\beta}}|r(\mathbf{kq})-\rho(\mathbf{T})|\nonumber\\
&\leq& C\frac{q_{1}q_{2}\log (T_{1}T_{2})}{T_{1}T_{2}}
      \sum_{\mathbf{0}\leq\mathbf{kq}\leq \mathbf{T}, \mathbf{kq}\neq \mathbf{0}\atop k_{1}q_{1}k_{2}q_{2}> (T_{1}T_{2})^{\beta}}|r(\mathbf{kq})-\rho(\mathbf{T})|\nonumber\\
&\leq&  C\frac{q_{1}q_{2}}{\beta T_{1}T_{2}}\sum_{\mathbf{0}\leq\mathbf{kq}\leq \mathbf{T}, \mathbf{kq}\neq \mathbf{0}\atop k_{1}q_{1}k_{2}q_{2}> (T_{1}T_{2})^{\beta}}|r(\mathbf{kq})\log (k_{1}q_{1}k_{2}q_{2})-r|\nonumber\\
&&    +Cr\frac{q_{1}q_{2}}{T_{1}T_{2}} \sum_{\mathbf{0}\leq\mathbf{kq}\leq \mathbf{T}, \mathbf{kq}\neq \mathbf{0}\atop k_{1}q_{1}k_{2}q_{2}> (T_{1}T_{2})^{\beta}}\left|1-\frac{\log (T_{1}T_{2})}{\log (k_{1}q_{1}k_{2}q_{2})}\right|.
\end{eqnarray}
By  Assumption  \textbf{A3}, the first term on the
right-hand-side of (\ref{eq409}) tends to 0 as $\mathbf{T}\to
\infty$. Furthermore, the second term of the right-hand-side of
(\ref{eq409}) also tends to 0 by an integral estimate as follows (
see also the proof of Lemma 6.4.1 of Leadbetter et al.\ (1983))
\begin{eqnarray*}
&&\frac{q_{1}q_{2}}{T_{1}T_{2}} \sum_{\mathbf{0}\leq\mathbf{kq}\leq \mathbf{T}, \mathbf{kq}\neq \mathbf{0}\atop k_{1}q_{1}k_{2}q_{2}> (T_{1}T_{2})^{\beta}}
     \left|1-\frac{\log (T_{1}T_{2})}{\log (k_{1}q_{1}k_{2}q_{2})}\right|\\
&&\leq\frac{r}{\log(T_{1}T_{2})^{\beta}} \frac{q_{1}q_{2}}{T_{1}T_{2}} \sum
\left| \log(k_{1}q_{1}k_{2}q_{2}) - \log (T_{1}T_{2})\right|\\
&&=\frac{r}{ \log(T_{1}T_{2})^{\beta}} \frac{q_{1}q_{2}}{T_{1}T_{2}} \sum
\left|\log\left(\frac{k_{1}q_{1}k_{2}q_{2}}{T_{1}T_{2}}\right)\right|\\
&&= O\left(\frac{r}{\log(T_{1}T_{2})^{\beta}}\int_0^1 \int_0^1 \log(xy)dxdy\right),
\end{eqnarray*}
which shows that $S_{\mathbf{T},12}\rightarrow 0$ as $\mathbf{T}\rightarrow\infty$. Thus, $S_{\mathbf{T},1}\rightarrow 0$ as $\mathbf{T}\rightarrow\infty$.

We consider the term $S_{\mathbf{T},2}$ and we will discuss it for two cases, the first for $(T_{1}T_{2})^{\beta}>T_{2}$, and the second for $(T_{1}T_{2})^{\beta}\leq T_{2}$.

For the case $(T_{1}T_{2})^{\beta}>T_{2}$, by the same arguments as for
(\ref{eq405}), we have
\begin{eqnarray*}
S_{\mathbf{T},2}
&=&C\frac{T_{1}}{q_{1}}\frac{T_{2}}{q_{2}}\sum_{0\leq k_{2}q_{2}\leq T_{2},k_{1}q_{1}=0}
|r(0,k_{2}q_{2})-\rho(\mathbf{T})|\exp\left(-\frac{u^{2}_{\mathbf{T}}}{1+\vartheta(T_{1}^{b},T_{2}^{b})}\right)\\
&\leq& C\frac{T_{1}}{q_{1}}\frac{T_{2}}{q_{2}}\sum_{0\leq k_{2}q_{2}\leq T_{2},k_{1}q_{1}=0}
     \exp\left(-\frac{u^{2}_{\mathbf{T}}}{1+\vartheta(\varepsilon,\varepsilon)}\right)\\
&\leq & C (T_{1}T_{2})^{\beta-\frac{1-\vartheta(\varepsilon,\varepsilon)}{1+\vartheta(\varepsilon,\varepsilon)}}(\log T_{1}T_{2})^{1/\alpha_{1}+2/\alpha_{2}}.
\end{eqnarray*}
Therefore, $S_{\mathbf{T},2}\rightarrow 0$ as $\mathbf{T}\rightarrow\infty$ in view of $\beta<\frac{1-\vartheta(\varepsilon,\varepsilon)}{1+\vartheta(\varepsilon,\varepsilon)}$.

For the second case $(T_{1}T_{2})^{\beta}\leq T_{2}$, split $S_{\mathbf{T},2}$ into two parts as
\begin{eqnarray*}
S_{\mathbf{T},2}= C\sum_{\mathbf{kq}\in \mathbf{O_{i}},\mathbf{lq}\in \mathbf{O_{j}}, \mathbf{kq}\neq \mathbf{lq}, \mathbf{1}\leq \mathbf{i}\neq \mathbf{j}\leq \mathbf{n} \atop 0< |k_{2}q_{2}-l_{2}q_{2}|\leq (T_{1}T_{2})^{\beta},k_{1}q_{1}=l_{1}q_{1}}
+ C\sum_{\mathbf{kq}\in \mathbf{O_{i}},\mathbf{lq}\in \mathbf{O_{j}},\mathbf{kq}\neq \mathbf{lq}, \mathbf{1}\leq \mathbf{i}\neq \mathbf{j}\leq \mathbf{n} \atop (T_{1}T_{2})^{\beta}< |k_{2}q_{2}-l_{2}q_{2}|\leq T_{2},k_{1}q_{1}=l_{1}q_{1}}=: S_{\mathbf{T},21}+S_{\mathbf{T},22}.
\end{eqnarray*}
For $S_{\mathbf{T},21}$, similarly to the derivation of
(\ref{eq405}) again, we have
\begin{eqnarray*}
S_{\mathbf{T},21}&\leq&C\frac{T_{1}}{q_{1}}\frac{T_{2}}{q_{2}}\sum_{0\leq k_{2}q_{2}\leq (T_{1}T_{2})^{\beta},k_{1}q_{1}=0}
|r(0,k_{2}q_{2})-\rho(\mathbf{T})|\exp\left(-\frac{u^{2}_{\mathbf{T}}}{1+\vartheta(T_{1}^{b},T_{2}^{b})}\right)\\
&\leq & C\frac{T_{1}}{q_{1}}\frac{T_{2}}{q_{2}}\sum_{0\leq k_{2}q_{2}\leq (T_{1}T_{2})^{\beta},k_{1}q_{1}=0}
     \exp\left(-\frac{u^{2}_{\mathbf{T}}}{1+\vartheta(\varepsilon,\varepsilon)}\right)\nonumber\\
&\leq & C (T_{1}T_{2})^{\beta-\frac{1-\vartheta(\varepsilon,\varepsilon)}{1+\vartheta(\varepsilon,\varepsilon)}}(\log T_{1}T_{2})^{1/\alpha_{1}+2/\alpha_{2}},
\end{eqnarray*}
which shows that $S_{\mathbf{T},21}\rightarrow 0$ as $\mathbf{T}\rightarrow\infty$.

For bound the term $S_{\mathbf{T},22}$, we need to define
$$\omega_{2}(\mathbf{t})=\max\{|r(0,t_{2})|, |\rho(\mathbf{T})|\}$$
and
$$
\theta_{2}(\mathbf{\mathbf{z}})=\sup_{\mathbf{0}\leq\mathbf{t}\leq \mathbf{T},\atop |t_{2}|>z_{1}z_{2}}\{\omega(\mathbf{t})\}.$$
By Assumption  \textbf{A3} again, we have also
$\theta_{2}(\mathbf{kq})\leq K/\log(T_{1}T_{2})^{\beta}$ and
$r(0,k_{2}q_{2})\log (T_{1}T_{2})\leq C$ for $k_{1}q_{1}=0$ and $k_{2}q_{2}>(T_{1}T_{2})^{\beta}$.
So by the same arguments as for (\ref{Tan2}), we have
$$\frac{(T_{1}T_{2})^{2}}{q_{1}^{2}q_{2}^{2}\log (T_{1}T_{2})}\exp\left(-\frac{u_{\mathbf{T}}^{2}}{1+\theta_{2}(T_{1}^{\beta},T_{2}^{\beta})}\right)=O(1)$$
and we  thus have
\begin{eqnarray*}
S_{\mathbf{T},22}&\leq&C\frac{T_{1}}{q_{1}}\frac{T_{2}}{q_{2}}\sum_{(T_{1}T_{2})^{\beta}< k_{2}q_{2}\leq T_{2},k_{1}q_{1}=0}
|r(0,k_{2}q_{2})-\rho(\mathbf{T})|\exp\left(-\frac{u^{2}_{\mathbf{T}}}{1+\theta_{2}(T_{1}^{\beta},T_{2}^{\beta})}\right)\\
&= & C\frac{(T_{1}T_{2})^{2}}{q_{1}^{2}q_{2}^{2}\log (T_{1}T_{2})}\exp\left(-\frac{u_{\mathbf{T}}^{2}}{1+\theta_{2}(T_{1}^{\beta},T_{2}^{\beta})}\right)\cdot \frac{q_{1}q_{2}\log (T_{1}T_{2})}{T_{1}T_{2}}
 \sum_{(T_{1}T_{2})^{\beta}< k_{2}q_{2}\leq T_{2},k_{1}q_{1}=0}|r(0,k_{2}q_{2})-\rho(\mathbf{T})|\\
&\leq &C \frac{q_{1}q_{2}\log (T_{1}T_{2})}{T_{1}T_{2}}
 \sum_{(T_{1}T_{2})^{\beta}< k_{2}q_{2}\leq T_{2},k_{1}q_{1}=0}|r(0,k_{2}q_{2})-\rho(\mathbf{T})|\\
&\leq &C \frac{q_{1}q_{2}\log (T_{1}T_{2})}{T_{1}T_{2}}
 \sum_{(T_{1}T_{2})^{\beta}< k_{2}q_{2}\leq T_{2},k_{1}q_{1}=0}(|r(0,k_{2}q_{2})|+\rho(\mathbf{T}))\\
 &\leq &C \frac{q_{1}q_{2}\log (T_{1}T_{2})}{T_{1}T_{2}}\frac{T_{2}}{q_{2}}\frac{1}{\log(T_{1}T_{2})}\\
 &=&C\frac{q_{1}}{T_{1}},
\end{eqnarray*}
which implies $S_{\mathbf{T},22}\rightarrow0$ as $\mathbf{T}\rightarrow\infty$.
Thus we have proved that  $S_{\mathbf{T},2}\rightarrow0$ as $\mathbf{T}\rightarrow\infty$.
By the same arguments, we can show that  $S_{\mathbf{T},3}\rightarrow0$ as $\mathbf{T}\rightarrow\infty$.
The proof of the lemma is complete.\hfill$\Box$\\

\textbf{Lemma B2}. {\sl Under the conditions of Lemma 3.3, we have
\begin{eqnarray}
\label{B2}
\sum_{\mathbf{\mathbf{kp}}\in \mathbf{O_{i}},\mathbf{\mathbf{lp}}\in \mathbf{O_{j}}\atop \mathbf{kp}\neq \mathbf{lp}, \mathbf{1}\leq \mathbf{i,j}\leq \mathbf{n} }|r(\mathbf{kp},\mathbf{lp})-\varrho(\mathbf{kp},\mathbf{lp})|
\int_{0}^{1}\frac{1}{\sqrt{1-r^{(h)}(\mathbf{kp},\mathbf{lp})}}\exp\left(-\frac{u_{\mathbf{T}}'^{2}}{1+r^{(h)}(\mathbf{kp},\mathbf{lp})}\right)dh
\rightarrow 0
\end{eqnarray}
as $\mathbf{T}\rightarrow\infty$.
}

\textbf{Proof:}
The proof is the same as that of Lemma B1, we omit the details.\\

\textbf{Lemma B3}. {\sl Under the conditions of Lemma 3.3, we have
\begin{eqnarray}
\label{eqB30}
\sum_{\mathbf{kq}\in \mathbf{O_{i}},\mathbf{lp}\in \mathbf{O_{j}}\atop \mathbf{kq}\neq \mathbf{lp}, \mathbf{1}\leq \mathbf{i,j}\leq \mathbf{n} }|r(\mathbf{kq},\mathbf{lp})-\varrho(\mathbf{kq},\mathbf{lp})|
\int_{0}^{1}\frac{1}{\sqrt{1-r^{(h)}(\mathbf{kq},\mathbf{lp})}}\exp\left(-\frac{u_{\mathbf{T}}'^{2}+u_{\mathbf{T}}^{2}}{2(1+r^{(h)}(\mathbf{kq},\mathbf{lp}))}\right)dh
\rightarrow 0
\end{eqnarray}
as $\mathbf{T}\rightarrow\infty$.
}

\textbf{Proof:} Recall that $\mathfrak{R}(p_{i})$, $i=1,2$ can be sparse grids or Pickands grids.
First, we consider the case that $\mathbf{kq},\mathbf{lp}$ in the
same interval $\mathbf{O_{i}}$.
Split the sum in (\ref{eqB30}) into two parts as
\begin{eqnarray}
\label{eqB31}
\sum_{\mathbf{kq},\mathbf{lp}\in \mathbf{O_{i}},\mathbf{kq}\neq\mathbf{lp}, \mathbf{i}=\mathbf{1},\cdots,\mathbf{n},\atop \max\{|l_{1}p_{1}-k_{1}q_{1}|,|l_{2}p_{2}-k_{2}q_{2}|\}\leq\varepsilon}
+\sum_{\mathbf{kq},\mathbf{lp}\in \mathbf{O_{i}},\mathbf{kq}\neq\mathbf{lp}, \mathbf{i}=\mathbf{1},\cdots,\mathbf{n}, \atop \max\{|l_{1}p_{1}-k_{1}q_{1}|,|l_{2}p_{2}-k_{2}q_{2}|\}>\varepsilon}=:W_{\mathbf{T},1}+W_{\mathbf{T},2}.
\end{eqnarray}
We deal with $W_{\mathbf{T},1}$.  For $\mathbf{kq},\mathbf{lp}$ in the
same interval $\mathbf{O_{i}}$, we  have
$\varrho(\mathbf{kq},\mathbf{lp})-r(\mathbf{kq},\mathbf{lq})=\rho(\mathbf{T})(1-r(\mathbf{kq},\mathbf{lq}))$.
By  Assumption {\bf A1}  we can also choose small enough $\varepsilon>0$ such that
$\varrho(\mathbf{kq},\mathbf{lp})=r(\mathbf{kq},\mathbf{lp})+(1-r(\mathbf{kq},\mathbf{lp}))\rho(\mathbf{T})\sim
r(\mathbf{kq},\mathbf{lp})$ for sufficiently large $\mathbf{T}$ and $\max\{|l_{1}p_{1}-k_{1}q_{1}|,|l_{2}p_{2}-k_{2}q_{2}|\}\leq\varepsilon$.
By the definitions of $u_{\mathbf{T}}$ and $u_{\mathbf{T}}'$, we have
\begin{eqnarray}
\label{eqB32}
v_{\mathbf{T}}^{2}:=\frac{1}{2}(u_{\mathbf{T}}^{2}+(u_{\mathbf{T}}')^{2})=2\log (T_{1}T_{2})-\log\log (T_{1}T_{2})+\log (p_{1}^{-1}p_{2}^{-1})+(1/\alpha_{1}+1/\alpha_{2})\log\log (T_{1}T_{2})+O(1).
\end{eqnarray}
Consequently, in view of (\ref{eqT20}), we obtain
\begin{eqnarray}
\label{eqB33}
W_{\mathbf{T},1}&\leq&C\sum_{\mathbf{kq},\mathbf{lp}\in \mathbf{O_{i}},\mathbf{kq}\neq\mathbf{lp}, \mathbf{i}=\mathbf{1},\cdots,\mathbf{n},\atop \max\{|l_{1}p_{1}-k_{1}q_{1}|,|l_{2}p_{2}-k_{2}q_{2}|\}\leq\varepsilon}
|r(\mathbf{kq},\mathbf{lp})-\varrho(\mathbf{kq},\mathbf{lp})|\frac{1}{\sqrt{1-r(\mathbf{kq},\mathbf{lp})}}
    \exp\left(-\frac{v_{\mathbf{T}}^{2}}{1+r(\mathbf{kq},\mathbf{lq})}\right)\nonumber\\
&\leq & C\rho(\mathbf{T})\sum_{\mathbf{kq},\mathbf{lp}\in \mathbf{O_{i}},\mathbf{kq}\neq\mathbf{lp}, \mathbf{i}=\mathbf{1},\cdots,\mathbf{n},\atop \max\{|l_{1}p_{1}-k_{1}q_{1}|,|l_{2}p_{2}-k_{2}q_{2}|\}\leq\varepsilon}\sqrt{1-r(\mathbf{kq},\mathbf{lp})}
\exp\left(-\frac{v^{2}_{\mathbf{T}}}{2}\right)\exp\left(-\frac{(1-r(\mathbf{kq},\mathbf{lp}))v^{2}_{\mathbf{T}}}{(1+r(\mathbf{kq},\mathbf{lp}))}\right)\nonumber\\
&\leq & C \rho(\mathbf{T})(T_{1}T_{2})^{-1}(p_{1}p_{2})^{1/2}(\log T_{1}T_{2})^{1/2-1/2\alpha_{1}-1/2\alpha_{2}}\times\nonumber\\
&&\sum_{\mathbf{kq},\mathbf{lp}\in \mathbf{O_{i}},\mathbf{kq}\neq\mathbf{lp}, \mathbf{i}=\mathbf{1},\cdots,\mathbf{n},\atop \max\{|l_{1}p_{1}-k_{1}q_{1}|,|l_{2}p_{2}-k_{2}q_{2}|\}\leq\varepsilon}\sqrt{1-r(\mathbf{kq}-\mathbf{lp})}
   \exp\left(-\frac{(1-r(\mathbf{kq}-\mathbf{lp}))u^{2}_{\mathbf{T}}}{2(1+r(\mathbf{kq}-\mathbf{lp}))}\right)\nonumber\\
&\leq & C(T_{1}T_{2})^{-1}(p_{1}p_{2})^{1/2}(\log T_{1}T_{2})^{-1/2-1/2\alpha_{1}-1/2\alpha_{2}}\times\nonumber\\
&&\sum_{\mathbf{kq},\mathbf{lp}\in \mathbf{O_{i}},\mathbf{kq}\neq\mathbf{lp}, \mathbf{i}=\mathbf{1},\cdots,\mathbf{n},\atop \max\{|l_{1}p_{1}-k_{1}q_{1}|,|l_{2}p_{2}-k_{2}q_{2}|\}\leq\varepsilon}\sqrt{|l_{1}p_{1}-k_{1}q_{1}|^{\alpha_{1}}+|l_{2}p_{2}-k_{2}q_{2}|^{\alpha_{2}}}\times\nonumber\\
 &&  \exp\left(-\frac{(|l_{1}p_{1}-k_{1}q_{1}|^{\alpha_{1}}+|l_{2}p_{2}-k_{2}q_{2}|^{\alpha_{2}})v^{2}_{\mathbf{T}}}{8}\right).
\end{eqnarray}
%Since $\mathfrak{R}(p_{i})$, $i=1,2$ are sparse grids, $p_{i}(\log
%T_{1}T_{2})^{1/\alpha_{i}}\rightarrow\infty$.
Noting that $q_{i}=\gamma_{i}\log(T_{1}T_{2})^{1/\alpha_{i}}$ and $\mathfrak{R}(p_{i})$, $i=1,2$ are sparse grids or Pickands grids, a direct calculation shows that
\begin{eqnarray*}
&&\sum_{\mathbf{kq},\mathbf{lp}\in \mathbf{O_{i}},\mathbf{kq}\neq\mathbf{lp}, \mathbf{i}=\mathbf{1},\cdots,\mathbf{n},\atop \max\{|l_{1}p_{1}-k_{1}q_{1}|,|l_{2}p_{2}-k_{2}q_{2}|\}\leq\varepsilon}\sqrt{|l_{1}p_{1}-k_{1}q_{1}|^{\alpha_{1}}+|l_{2}p_{2}-k_{2}q_{2}|^{\alpha_{2}}}\times\nonumber\\
 &&  \ \ \ \ \ \ \ \ \ \ \  \exp\left(-\frac{(|l_{1}p_{1}-k_{1}q_{1}|^{\alpha_{1}}+|l_{2}p_{2}-k_{2}q_{2}|^{\alpha_{2}})v^{2}_{\mathbf{T}}}{8}\right)\\
 &&\ \ \ \ \ \leq CT_{1}T_{2}p_{1}^{-1}p_{2}^{-1}\sum_{0<k_{1}q_{1}<\varepsilon, 0<k_{2}q_{2}<\varepsilon}
\exp\left(-\frac{1}{4}[(k_{1}q_{1})^{\alpha_{1}}+(k_{2}q_{2})^{\alpha_{2}}]\log(T_{1}T_{2})\right)\\
&&\ \ \ \ \ \leq CT_{1}T_{2}p_{1}^{-1}p_{2}^{-1}\sum_{k_{1}=1}^{\infty}e^{-\frac{1}{4}(k_{1}\gamma_{1})^{\alpha_{1}}}
\sum_{k_{2}=1}^{\infty}e^{-\frac{1}{4}(k_{2}\gamma_{2})^{\alpha_{2}}}\\
&&\ \ \ \ \ \leq CT_{1}T_{2}p_{1}^{-1}p_{2}^{-1},
\end{eqnarray*}
which combine with (\ref{eqB33}) shows that $W_{\mathbf{T},1}\rightarrow0$ as $\mathbf{T}\rightarrow\infty$. \\
Using the fact that $v_{\mathbf{T}}\thicksim (2\log
T_{1}T_{2})^{1/2}$, we obtain
\begin{eqnarray}
\label{eqB34}
W_{\mathbf{T},2}&\leq&C\sum_{\mathbf{kq},\mathbf{lp}\in \mathbf{O_{i}},\mathbf{kq}\neq\mathbf{lp}, \mathbf{i}=\mathbf{1},\cdots,\mathbf{n}, \atop \max\{|l_{1}p_{1}-k_{1}q_{1}|,|l_{2}p_{2}-k_{2}q_{2}|\}>\varepsilon}
      |r(\mathbf{kq},\mathbf{lp})-\varrho(\mathbf{kq},\mathbf{lp})|\exp\left(-\frac{v^{2}_{\mathbf{T}}}{1+\varpi(\mathbf{kq},\mathbf{lp})}\right)\nonumber\\
&\leq & C\sum_{\mathbf{kq},\mathbf{lp}\in \mathbf{O_{i}},\mathbf{kq}\neq\mathbf{lp},  \mathbf{i}=\mathbf{1},\cdots,\mathbf{n},\atop \max\{|l_{1}p_{1}-k_{1}q_{1}|,|l_{2}p_{2}-k_{2}q_{2}|\}>\varepsilon}
     \exp\left(-\frac{v^{2}_{\mathbf{T}}}{1+\vartheta(\varepsilon,\varepsilon)}\right)\nonumber\\
&\leq & C\frac{T_{1}}{p_{1}}\frac{T_{2}}{p_{2}}\exp\left(-\frac{u^{2}_{\mathbf{T}}}{1+\vartheta(\varepsilon,\varepsilon)}\right)
     \sum_{0\leq k_{1}q_{1}\leq T_{1}^{a}, 0\leq k_{2}q_{2}\leq T_{2}^{a}}1\nonumber\\
&\leq & C \frac{T_{1}}{p_{1}}\frac{T_{2}}{p_{2}}(T_{1}T_{2})^{-\frac{2}{1+\vartheta(\varepsilon,\varepsilon)}}
    \sum_{0\leq k_{1}q_{1}\leq T_{1}^{a}, 0\leq k_{2}q_{2}\leq T_{2}^{a}}1\nonumber\\
&\leq & C (T_{1}T_{2})^{a-\frac{1-\vartheta(\varepsilon,\varepsilon)}{1+\vartheta(\varepsilon,\varepsilon)}}(\log T_{1}T_{2})^{1/\alpha_{1}+1/\alpha_{2}}(p_{1}p_{2})^{-1}.
\end{eqnarray}
Thus, $W_{\mathbf{T},2}\rightarrow0$ as $\mathbf{T}\rightarrow\infty$ by virtue of $a<\frac{1-\vartheta(\varepsilon,\varepsilon)}{1+\vartheta(\varepsilon,\varepsilon)}$ again.

Second, we deal with the case that $\mathbf{kq}\in \mathbf{O_{i}}$
and $\mathbf{lp}\in \mathbf{O_{j}}$, $\mathbf{i}\neq \mathbf{j}$.
Note that in this case, the distance between the points in any two
rectangles $\mathbf{O_{i}}$ and $\mathbf{O_{j}}$ is large than
$T_{1}^{b}$ or $T_{2}^{b}$ and $\varrho(\mathbf{kq},\mathbf{lp})=\rho(\mathbf{T})$ for  $\mathbf{kq}\in \mathbf{O_{i}}$
and $\mathbf{lp}\in \mathbf{O_{j}}$, $\mathbf{i}\neq \mathbf{j}$.

Obviously, the sum in (\ref{eqB30}) is at most
\begin{eqnarray}
\label{eqB35}
C\sum_{\mathbf{kq}\in \mathbf{O_{i}},\mathbf{lp}\in \mathbf{O_{j}}\atop \mathbf{kq}\neq \mathbf{lp}, \mathbf{1}\leq \mathbf{i}\neq \mathbf{j}\leq \mathbf{n} }
|r(\mathbf{kq},\mathbf{lp})-\rho(\mathbf{T})|\exp\left(-\frac{v^{2}_{\mathbf{T}}}{1+\varpi(\mathbf{kq},\mathbf{lp})}\right).
%&\leq&C\sum_{\mathbf{kq}\in \mathbf{O_{i}},\mathbf{lq}\in \mathbf{O_{j}}\atop \mathbf{kq}\neq \mathbf{lq}, \mathbf{1}\leq \mathbf{i}\neq \mathbf{j}\leq \mathbf{n} }
%|r(\mathbf{kq},\mathbf{lq})-\rho(\mathbf{T})|\exp\left(-\frac{u^{2}_{\mathbf{T}}}{1+\vartheta(T_{1}^{b},T_{2}^{b})}\right)\nonumber\\
%&\leq&  C\frac{T_{1}}{q_{1}}\frac{T_{2}}{q_{2}}\sum_{\mathbf{0}\leq\mathbf{kq}\leq \mathbf{T}, \mathbf{kq}\neq \mathbf{0}}
%|r(\mathbf{kq})-\rho(\mathbf{T})|\exp\left(-\frac{u^{2}_{\mathbf{T}}}{1+\vartheta(T_{1}^{b},T_{2}^{b})}\right)
\end{eqnarray}
Split the sum of (\ref{eqB35}) into three
parts, the first for $|k_{1}q_{1}-l_{1}p_{1}|>0$ and $|k_{2}q_{2}-l_{2}p_{2}|>0$, the second for
$k_{1}q_{1}-l_{1}p_{1}=0$ and $|k_{2}q_{2}-l_{2}p_{2}|>0$, the third for $k_{2}q_{2}-l_{2}p_{2}=0$ and $|k_{1}q_{1}-l_{1}p_{1}|>0$ and
denote them by $M_{\mathbf{T},i}$, $i=1,2,3,$ respectively. Let $\beta$ be chosen as before, ie,
$0<b<a<\beta<\frac{1-\vartheta(\varepsilon,\varepsilon)}{1+\vartheta(\varepsilon,\varepsilon)}$
for all sufficiently large $\mathbf{T}$.

We consider the term $M_{\mathbf{T},1}$ and
split it into two parts as
\begin{eqnarray*}
M_{\mathbf{T},1}= C\sum_{\mathbf{kq}\in \mathbf{O_{i}},\mathbf{lp}\in \mathbf{O_{j}}, \mathbf{kq}\neq \mathbf{lp}, \mathbf{1}\leq \mathbf{i}\neq \mathbf{j}\leq \mathbf{n}\atop |k_{1}q_{1}-l_{1}q_{1}||k_{2}q_{2}-l_{2}q_{2}|\leq (T_{1}T_{2})^{\beta}}
+ C\sum_{\mathbf{kq}\in \mathbf{O_{i}},\mathbf{lp}\in \mathbf{O_{j}}, \mathbf{kq}\neq \mathbf{lp}, \mathbf{1}\leq \mathbf{i}\neq \mathbf{j}\leq \mathbf{n}\atop |k_{1}q_{1}-l_{1}q_{1}||k_{2}q_{2}-l_{2}q_{2}|> (T_{1}T_{2})^{\beta}}=: M_{\mathbf{T},11}+M_{\mathbf{T},12}.
\end{eqnarray*}
For $M_{\mathbf{T},11}$, with the similar derivation as for
(\ref{eqB34}), we have
\begin{eqnarray}
\label{eqB36}
M_{\mathbf{T},11}
&\leq & C\frac{T_{1}}{p_{1}}\frac{T_{2}}{p_{2}}\sum_{0\leq k_{1}q_{1}\leq T_{1}, 0\leq k_{2}q_{2}\leq T_{2},\atop k_{1}q_{1}k_{2}q_{2}\leq (T_{1}T_{2})^{\beta}}
     \exp\left(-\frac{v^{2}_{\mathbf{T}}}{1+\vartheta(\varepsilon,\varepsilon)}\right)\nonumber\\
&\leq & C (T_{1}T_{2})^{\beta-\frac{1-\vartheta(\varepsilon,\varepsilon)}{1+\vartheta(\varepsilon,\varepsilon)}}(\log T_{1}T_{2})^{1/\alpha_{1}+1/\alpha_{2}}(p_{1}p_{2})^{-1}.
\end{eqnarray}
Thus, $M_{\mathbf{T},11}\rightarrow0$ as $\mathbf{T}\rightarrow\infty$, since $\beta<\frac{1-\vartheta(\varepsilon,\varepsilon)}{1+\vartheta(\varepsilon,\varepsilon)}$.\\
For $M_{\mathbf{T},12}$, we need also more precise estimation. Recall that
$$\omega_{1}(\mathbf{t})=\max\{|r(\mathbf{t})|, |\rho(\mathbf{T})|\}\ \ \mbox{and}\ \
\theta_{1}(\mathbf{\mathbf{z}})=\sup_{\mathbf{0}\leq\mathbf{t}\leq \mathbf{T},\atop |t_{1}t_{2}|>z_{1}z_{2}}\{\omega_{1}(\mathbf{t})\}.$$
Now using (\ref{eqB32}) again, by the same arguments as for (\ref{Tan2}), we
obtain
\begin{eqnarray}
\label{eqB37}
&&\frac{(T_{1}T_{2})^{2}}{q_{1}q_{2}p_{1}p_{2}\log (T_{1}T_{2})}\exp\left(-\frac{v_{\mathbf{T}}^{2}}{1+\theta_{1}(T_{1}^{\beta},T_{2}^{\beta})}\right)\nonumber\\
&&\leq \frac{(T_{1}T_{2})^{2}}{q_{1}q_{2}p_{1}p_{2}\log (T_{1}T_{2})}\exp\left(-\frac{v_{\mathbf{T}}^{2}}{1+K/\log(T_{1}T_{2})^{\beta}}\right)\nonumber\\
&&\thicksim \frac{(T_{1}T_{2})^{2}}{q_{1}q_{2}p_{1}p_{2}\log (T_{1}T_{2})}\left((T_{1}T_{2})^{-2}(\log T_{1}T_{2}) (\log T_{1}T_{2})^{-(1/\alpha_{1}+1/\alpha_{2})}(p_{1}p_{2})^{-1}\right)^{\frac{1}{1+K/\log(T_{1}T_{2})^{\beta}}}\nonumber\\
&&=O(1)
\end{eqnarray}
and then we thus have
\begin{eqnarray}
\label{eqB38}
M_{\mathbf{T},12}&\leq&C\sum_{\mathbf{kq}\in \mathbf{O_{i}},\mathbf{lp}\in \mathbf{O_{j}}, \mathbf{kq}\neq \mathbf{lp}, \mathbf{1}\leq \mathbf{i}\neq \mathbf{j}\leq \mathbf{n}\atop |k_{1}q_{1}-l_{1}q_{1}||k_{2}q_{2}-l_{2}q_{2}|> (T_{1}T_{2})^{\beta} }
|r(\mathbf{kq}-\mathbf{lp})-\rho(\mathbf{T})|\exp\left(-\frac{v^{2}_{\mathbf{T}}}{1+\theta_{1}(T_{1}^{\beta},T_{2}^{\beta})}\right)\nonumber\\
&=& C\frac{(T_{1}T_{2})^{2}}{q_{1}q_{2}p_{1}p_{2}\log (T_{1}T_{2})}\exp\left(-\frac{v_{\mathbf{T}}^{2}}{1+\theta_{1}(T_{1}^{\beta},T_{2}^{\beta})}\right)\times\nonumber\\
 &&\ \ \ \frac{q_{1}q_{2}p_{1}p_{2}\log (T_{1}T_{2})}{(T_{1}T_{2})^{2}}
     \sum_{\mathbf{kq}\in \mathbf{O_{i}},\mathbf{lp}\in \mathbf{O_{j}}, \mathbf{kq}\neq \mathbf{lp}, \mathbf{1}\leq \mathbf{i}\neq \mathbf{j}\leq \mathbf{n}\atop |k_{1}q_{1}-l_{1}q_{1}||k_{2}q_{2}-l_{2}q_{2}|> (T_{1}T_{2})^{\beta} }|r(\mathbf{kq}-\mathbf{lp})-\rho(\mathbf{T})|\nonumber\\
&\leq& C\frac{q_{1}q_{2}p_{1}p_{2}\log (T_{1}T_{2})}{(T_{1}T_{2})^{2}}
     \sum_{\mathbf{kq}\in \mathbf{O_{i}},\mathbf{lp}\in \mathbf{O_{j}}, \mathbf{kq}\neq \mathbf{lp}, \mathbf{1}\leq \mathbf{i}\neq \mathbf{j}\leq \mathbf{n}\atop |k_{1}q_{1}-l_{1}q_{1}||k_{2}q_{2}-l_{2}q_{2}|> (T_{1}T_{2})^{\beta} }|r(\mathbf{kq}-\mathbf{lp})-\rho(\mathbf{T})|\nonumber\\
&\leq&  C\frac{q_{1}q_{2}p_{1}p_{2}}{\beta (T_{1}T_{2})^{2}} \sum_{\mathbf{kq}\in \mathbf{O_{i}},\mathbf{lp}\in \mathbf{O_{j}}, \mathbf{kq}\neq \mathbf{lp}, \mathbf{1}\leq \mathbf{i}\neq \mathbf{j}\leq \mathbf{n}\atop |k_{1}q_{1}-l_{1}q_{1}||k_{2}q_{2}-l_{2}q_{2}|> (T_{1}T_{2})^{\beta} }|r(\mathbf{kq}-\mathbf{lp})\log ((k_{1}q_{1}-l_{1}p_{1})(k_{2}q_{2}-l_{2}p_{2}))-r|\nonumber\\
&&    +Cr\frac{q_{1}q_{2}p_{1}p_{2}}{ (T_{1}T_{2})^{2}}  \sum_{\mathbf{kq}\in \mathbf{O_{i}},\mathbf{lp}\in \mathbf{O_{j}}, \mathbf{kq}\neq \mathbf{lp}, \mathbf{1}\leq \mathbf{i}\neq \mathbf{j}\leq \mathbf{n}\atop |k_{1}q_{1}-l_{1}q_{1}||k_{2}q_{2}-l_{2}q_{2}|> (T_{1}T_{2})^{\beta} }\left|1-\frac{\log (T_{1}T_{2})}{\log ((k_{1}q_{1}-l_{1}p_{1})(k_{2}q_{2}-l_{2}p_{2}))}\right|.
\end{eqnarray}
By  Assumption  \textbf{A3}, the first term on the
right-hand-side of (\ref{eqB38}) tends to 0 as $\mathbf{T}\to
\infty$. Furthermore, the second term of the right-hand-side of
(\ref{eqB38}) also tends to 0 by an integral estimate as for Lemma B1. Thus $M_{\mathbf{T},12}\rightarrow 0$ as $\mathbf{T}\rightarrow\infty$ and then
$M_{\mathbf{T},1}\rightarrow 0$ as $\mathbf{T}\rightarrow\infty$.

We consider the term $M_{\mathbf{T},2}$ now. As the the proof of the previous lemma,  we also discuss it for two cases, the first for $(T_{1}T_{2})^{\beta}>T_{2}$, and the second for $(T_{1}T_{2})^{\beta}\leq T_{2}$.

For the case $(T_{1}T_{2})^{\beta}>T_{2}$, by the same arguments as for
(\ref{eqB34}), we have
\begin{eqnarray*}
M_{\mathbf{T},2}
&=&C\sum_{\mathbf{kq}\in \mathbf{O_{i}},\mathbf{lp}\in \mathbf{O_{j}}, \mathbf{kq}\neq \mathbf{lp}, \mathbf{1}\leq \mathbf{i}\neq \mathbf{j}\leq \mathbf{n}\atop k_{1}q_{1}-l_{1}p_{1}=0}
|r(0,k_{2}q_{2}-l_{2}p_{2})-\rho(\mathbf{T})|\exp\left(-\frac{v^{2}_{\mathbf{T}}}{1+\vartheta(T_{1}^{b},T_{2}^{b})}\right)\\
&\leq& C\exp\left(-\frac{v^{2}_{\mathbf{T}}}{1+\vartheta(\varepsilon,\varepsilon)}\right)
\sum_{\mathbf{kq}\in \mathbf{O_{i}},\mathbf{lp}\in \mathbf{O_{j}}, \mathbf{kq}\neq \mathbf{lp}, \mathbf{1}\leq \mathbf{i}\neq \mathbf{j}\leq \mathbf{n}\atop k_{1}q_{1}-l_{1}p_{1}=0}|r(0,k_{2}q_{2}-l_{2}p_{2})-\rho(\mathbf{T})|\\
&\leq & C (T_{1}T_{2})^{\beta-\frac{1-\vartheta(\varepsilon,\varepsilon)}{1+\vartheta(\varepsilon,\varepsilon)}}(\log T_{1}T_{2})^{1/\alpha_{2}}(p_{1}p_{2})^{-1},
\end{eqnarray*}
which shows that $M_{\mathbf{T},2}\rightarrow 0$ as $\mathbf{T}\rightarrow\infty$, where in the last step, we use the fact that the number of $(k_{1},l_{1})$ such that $k_{1}q_{1}-l_{1}p_{1}=0$ does not exceed $T_{1}/p_{1}$ and $(T_{1}T_{2})^{\beta}>T_{2}$.

For the second case $(T_{1}T_{2})^{\beta}\leq T_{2}$, split $M_{\mathbf{T},2}$ into two parts as
\begin{eqnarray*}
M_{\mathbf{T},2}= C\sum_{\mathbf{kq}\in \mathbf{O_{i}},\mathbf{lp}\in \mathbf{O_{j}}, \mathbf{kq}\neq \mathbf{lp}, \mathbf{1}\leq \mathbf{i}\neq \mathbf{j}\leq \mathbf{n}\atop0< |k_{2}q_{2}-l_{2}p_{2}|\leq (T_{1}T_{2})^{\beta},k_{1}q_{1}=l_{1}p_{1}}
+ C\sum_{\mathbf{kq}\in \mathbf{O_{i}},\mathbf{lp}\in \mathbf{O_{j}}, \mathbf{kq}\neq \mathbf{lp}, \mathbf{1}\leq \mathbf{i}\neq \mathbf{j}\leq \mathbf{n}\atop(T_{1}T_{2})^{\beta}< |k_{2}q_{2}-l_{2}p_{2}|\leq T_{2},k_{1}q_{1}=l_{1}p_{1}}=: M_{\mathbf{T},21}+M_{\mathbf{T},22}.
\end{eqnarray*}
For $M_{\mathbf{T},21}$, similarly to the derivation of
(\ref{eqB34}) again, we have
\begin{eqnarray*}
M_{\mathbf{T},21}&\leq&C\sum_{\mathbf{kq}\in \mathbf{O_{i}},\mathbf{lp}\in \mathbf{O_{j}}, \mathbf{kq}\neq \mathbf{lp}, \mathbf{1}\leq \mathbf{i}\neq \mathbf{j}\leq \mathbf{n}\atop 0< |k_{2}q_{2}-l_{2}p_{2}|\leq (T_{1}T_{2})^{\beta}, k_{1}q_{1}-l_{1}p_{1}=0}
|r(0,k_{2}q_{2}-l_{2}p_{2})-\rho(\mathbf{T})|\exp\left(-\frac{v^{2}_{\mathbf{T}}}{1+\vartheta(T_{1}^{b},T_{2}^{b})}\right)\\
&\leq& C\exp\left(-\frac{v^{2}_{\mathbf{T}}}{1+\vartheta(\varepsilon,\varepsilon)}\right)
\sum_{\mathbf{kq}\in \mathbf{O_{i}},\mathbf{lp}\in \mathbf{O_{j}}, \mathbf{kq}\neq \mathbf{lp}, \mathbf{1}\leq \mathbf{i}\neq \mathbf{j}\leq \mathbf{n}\atop 0< |k_{2}q_{2}-l_{2}p_{2}|\leq (T_{1}T_{2})^{\beta}, k_{1}q_{1}-l_{1}p_{1}=0}|r(0,k_{2}q_{2}-l_{2}p_{2})-\rho(\mathbf{T})|\\
&\leq & C (T_{1}T_{2})^{\beta-\frac{1-\vartheta(\varepsilon,\varepsilon)}{1+\vartheta(\varepsilon,\varepsilon)}}(\log T_{1}T_{2})^{1/\alpha_{2}}(p_{1}p_{2})^{-1},
\end{eqnarray*}
which shows that $M_{\mathbf{T},21}\rightarrow 0$ as $\mathbf{T}\rightarrow\infty$.

For $M_{\mathbf{T},22}$, we recall that
$$\omega_{2}(\mathbf{t})=\max\{|r(0,t_{2})|, |\rho(\mathbf{T})|\}\ \
\mbox{and}\ \
\theta_{2}(\mathbf{\mathbf{z}})=\sup_{\mathbf{0}\leq\mathbf{t}\leq \mathbf{T},\atop |t_{2}|>z_{1}z_{2}}\{\omega(\mathbf{t})\}.$$
So by the same arguments as for (\ref{Tan2}), we have
$$\frac{(T_{1}T_{2})^{2}}{q_{1}q_{2}p_{1}p_{2}\log (T_{1}T_{2})}\exp\left(-\frac{v_{\mathbf{T}}^{2}}{1+\theta_{2}(T_{1}^{\beta},T_{2}^{\beta})}\right)=O(1)$$
and we  thus have
\begin{eqnarray*}
M_{\mathbf{T},22}&\leq&C
\sum_{\mathbf{kq}\in \mathbf{O_{i}},\mathbf{lp}\in \mathbf{O_{j}}, \mathbf{kq}\neq \mathbf{lp}, \mathbf{1}\leq \mathbf{i}\neq \mathbf{j}\leq \mathbf{n}\atop (T_{1}T_{2})^{\beta}<|k_{2}q_{2}-l_{2}p_{2}| , k_{1}q_{1}-l_{1}p_{1}=0}
|r(0,k_{2}q_{2}-l_{2}p_{2})-\rho(\mathbf{T})|\exp\left(-\frac{v^{2}_{\mathbf{T}}}{1+\theta_{2}(T_{1}^{\beta},T_{2}^{\beta})}\right)\\
&= & C\frac{(T_{1}T_{2})^{2}}{q_{1}q_{2}p_{1}p_{2}\log (T_{1}T_{2})}\exp\left(-\frac{v_{\mathbf{T}}^{2}}{1+\theta_{2}(T_{1}^{\beta},T_{2}^{\beta})}\right)\times\\
&&  \frac{q_{1}q_{2}p_{1}p_{2}\log (T_{1}T_{2})}{(T_{1}T_{2})^{2}}
 \sum_{\mathbf{kq}\in \mathbf{O_{i}},\mathbf{lp}\in \mathbf{O_{j}}, \mathbf{kq}\neq \mathbf{lp}, \mathbf{1}\leq \mathbf{i}\neq \mathbf{j}\leq \mathbf{n}\atop  (T_{1}T_{2})^{\beta}<|k_{2}q_{2}-l_{2}p_{2}| , k_{1}q_{1}-l_{1}p_{1}=0}|r(0,k_{2}q_{2}-l_{2}p_{2})-\rho(\mathbf{T})|\\
&\leq &C \frac{q_{1}q_{2}p_{1}p_{2}\log (T_{1}T_{2})}{(T_{1}T_{2})^{2}}
 \sum_{\mathbf{kq}\in \mathbf{O_{i}},\mathbf{lp}\in \mathbf{O_{j}}, \mathbf{kq}\neq \mathbf{lp}, \mathbf{1}\leq \mathbf{i}\neq \mathbf{j}\leq \mathbf{n}\atop  (T_{1}T_{2})^{\beta}<|k_{2}q_{2}-l_{2}p_{2}| , k_{1}q_{1}-l_{1}p_{1}=0}|r(0,k_{2}q_{2}-l_{2}p_{2})-\rho(\mathbf{T})|\\
&\leq &C \frac{q_{1}q_{2}p_{1}p_{2}\log (T_{1}T_{2})}{(T_{1}T_{2})^{2}}
 \sum_{\mathbf{kq}\in \mathbf{O_{i}},\mathbf{lp}\in \mathbf{O_{j}}, \mathbf{kq}\neq \mathbf{lp}, \mathbf{1}\leq \mathbf{i}\neq \mathbf{j}\leq \mathbf{n}\atop  (T_{1}T_{2})^{\beta}<|k_{2}q_{2}-l_{2}p_{2}| , k_{1}q_{1}-l_{1}p_{1}=0}(|r(0,k_{2}q_{2}-l_{2}p_{2})|+\rho(\mathbf{T}))\\
 &\leq &C \frac{q_{1}q_{2}p_{1}p_{2}\log (T_{1}T_{2})}{(T_{1}T_{2})^{2}}\frac{T_{2}^{2}T_{1}}{p_{1}q_{2}p_{2}}\frac{1}{\log(T_{1}T_{2})}\\
 &=&C\frac{q_{1}}{T_{1}},
\end{eqnarray*}
which implies $M_{\mathbf{T},22}\rightarrow0$ as $\mathbf{T}\rightarrow\infty$.
Now we have showed that  $M_{\mathbf{T},2}\rightarrow0$ as $\mathbf{T}\rightarrow\infty$.
By the same arguments, we can show that  $M_{\mathbf{T},3}\rightarrow0$ as $\mathbf{T}\rightarrow\infty$.
The proof of the lemma is complete.\hfill$\Box$

\bigskip
{\bf Acknowledgement}:
We would like to thank the referee for several important suggestions
and corrections.

\end{document}